\documentclass[11pt,a4paper]{amsart}

\date{15 March  2006. Revised 21 January 2008}
\usepackage{amsfonts,amsmath,amssymb,indentfirst,mathrsfs,amscd}
\usepackage[mathscr]{eucal}
\usepackage{color,array}
\usepackage{colortbl}  

\DeclareMathOperator{\Quot}{Quot\,}        
\DeclareMathOperator{\img}{Im\,}
\DeclareMathOperator{\Gr}{Gr}        

\DeclareMathAlphabet{\mathpzc}{OT1}{pzc}{m}{it}
\DeclareMathOperator{\rk}{rk\,}
\DeclareMathOperator{\im}{im\,}
\DeclareMathOperator{\Tor}{Tor\,}
\DeclareMathOperator{\ParAut}{ParAut\,}
\DeclareMathOperator{\coker}{coker\,}

\DeclareMathOperator{\pdeg}{pardeg\,}
\DeclareMathOperator{\pmu}{par \mu\,}
\DeclareMathOperator{\PH}{ParHom\,}
\DeclareMathOperator{\SPH}{SParHom\,}
\DeclareMathOperator{\PE}{ParEnd\,}
\DeclareMathOperator{\SPE}{SParEnd\,}
\DeclareMathOperator{\ad}{ad}
\DeclareMathOperator{\orb}{orb}
\DeclareMathOperator{\kod}{kod}

\DeclareMathOperator{\Hom}{Hom\,}

\DeclareMathOperator{\id}{Id\,}
\DeclareMathOperator{\tr}{Tr\,}

\DeclareMathOperator{\Ext}{Ext}

\DeclareMathOperator{\U}{U}
\DeclareMathOperator{\SU}{SU}

\DeclareMathOperator{\GL}{GL}
\DeclareMathOperator{\length}{length}
\newcommand{\lie}{\mathfrak}
\newtheorem{thm}{Theorem}[section]
\newtheorem{prop}[thm]{Proposition}
\newtheorem{lem}[thm]{Lemma}
\newtheorem{cor}[thm]{Corollary}
\theoremstyle{definition}
\newtheorem{defn}[thm]{Definition}

\newtheorem{assumption}[thm]{Assumption}
\theoremstyle{remark}
\newtheorem{rmk}[thm]{Remark}

\newcommand{\lto}{\longrightarrow}
\newcommand{\iacute}{\'{\i}}

\newcommand{\imat}{\sqrt{-1}}
\newcommand{\dbar}{\bar{\partial}}
\newcommand{\cC}{\mathcal{C}}

\newcommand{\cE}{\mathcal{E}}

\newcommand{\cO}{\mathcal{O}}
\newcommand{\cF}{\mathcal{F}}
\newcommand{\cH}{\mathcal{H}}
\newcommand{\cM}{\mathcal{M}}
\newcommand{\cN}{\mathcal{N}}

\newcommand{\cQ}{\mathcal{Q}}
\newcommand{\cR}{\mathcal{R}}
\newcommand{\cS}{\mathcal{S}}
\newcommand{\cT}{\mathcal{T}}
\newcommand{\cU}{\mathcal{U}}

\newcommand{\CC}{\mathbb{C}}

\newcommand{\HH}{\mathbb{H}}
\newcommand{\PP}{\mathbb{P}}
\newcommand{\RR}{\mathbb{R}}
\newcommand{\ZZ}{\mathbb{Z}}
\newcommand{\AAA}{{\curly A}}
\newcommand{\CCC}{{\curly C}}

\newcommand{\GGG}{{\curly G}}
\newcommand{\HHH}{{\curly H}}

\DeclareFontFamily{OT1}{rsfs}{}
\DeclareFontShape{OT1}{rsfs}{n}{it}{<->rsfs10}{}
\DeclareMathAlphabet{\curly}{OT1}{rsfs}{n}{it}

\newcommand{\abs}[1]{\lvert #1\rvert}

\title{Moduli spaces of parabolic $\U(p,q)$-Higgs bundles}

\author{O. Garc\'{\i}a-Prada}
  \address{Instituto de Ciencias Matem{\'a}ticas CSIC-UAM-UCM-UC3M \\
  Consejo Superior de Investigaciones Cient{\'\i}ficas  \\ Serrano 121
  \\ 28006 Madrid \\ Spain}
   \email{oscar.garcia-prada@uam.es}
\author{M. Logares}
   \address{Departamento de Matematica Pura\\
Facultade de Ciencias\\
Rua do Campo Alegre 687\\
4169-007 Porto\\
Portugal }
   \email{mlogares@fc.up.pt}
\author{Vicente Mu\~noz}
  \address{Instituto de Ciencias Matem{\'a}ticas CSIC-UAM-UCM-UC3M \\
  Consejo Superior de Investigaciones Cient{\'\i}ficas \\ Serrano 113 bis
  \\ 28006 Madrid \\ Spain}
  \address{Facultad de Matem\'{a}ticas \\ Universidad Complutense
  de Madrid \\ Plaza Ciencias 3
  \\ 28040 Madrid \\ Spain}
  \email{vicente.munoz@imaff.cfmac.csic.es}

\thanks{Partially supported by Ministerio de Educaci{\'o}n
y Ciencia (Spain) through Project MTM2004-07090-C03-01.}



\begin{document}

\begin{abstract}
Using the $L^2$-norm of the Higgs field as a Morse function, we
study the moduli space of
parabolic $\U(p,q)$-Higgs bundles over a Riemann surface with a
finite number of marked points, under certain genericity
conditions on the parabolic structure. When the parabolic degree is zero this space is homeomorphic
to the moduli space of representations of the fundamental group of
the punctured  surface in $\U(p,q)$, with fixed compact holonomy
classes around the marked points. By means of this homeomorphism we count the number of connected components of this moduli space of representations. Finally, we apply our results to the
study of representations of the fundamental group of elliptic
surfaces of general type.
\end{abstract}

\maketitle

\section{Introduction}\label{sec:intro}

A parabolic vector bundle over a  compact Riemann surface with
marked points consists of a vector bundle, equipped  with a
weighted flag structure on  the  fibre over each marked  point.
These objects were introduced by Seshadri \cite{se} in relation to
certain desingularisations of the moduli space of semistable
vector bundles. It turns out that, similarly to the Narasimhan and
Seshadri correspondence \cite{nas,d} between stable vector bundles
and representations of the fundamental group of the surface in the
unitary group $\U(n)$, there is an analogous correspondence,
proved by Metha and Seshadri \cite{ms} (see also \cite{b}),
relating  stable parabolic
bundles to unitary representations of the fundamental group of the
punctured surface with a fixed holonomy class around each
marked point.

In order to study representations of the fundamental group of the
punctured surface in $\GL(n,\CC)$ one has to  consider  parabolic
Higgs bundles. These are pairs consisting of a parabolic  vector
bundle and  a meromorphic endomorphism valued one-form with a simple
pole along each  marked point, whose  residue is  nilpotent
with respect to the flag.
Moduli spaces of parabolic Higgs bundles provide interesting examples of
hyperk\"ahler manifolds.  This theory, studied by Simpson in
\cite{s2} and others \cite{by,k,na,nst}, generalizes the
non-parabolic Higgs bundle theory studied by Hitchin \cite{h},
Donaldson \cite{d2}, Simpson \cite{s1} and Corlette \cite{c}.

In this paper we study parabolic $\U(p,q)$-Higgs bundles. These
are the objects that correspond to representations of the
fundamental group of the punctured surface in $\U(p,q)$,  with
fixed  compact holonomy classes around the marked points. Our approach
combines the techniques used in \cite{bgg} in the study of
$\U(p,q)$-Higgs bundles in the non-parabolic case as well as those
used in \cite{ggm} to study the topology of moduli spaces of
$\GL(n,\CC)$-parabolic Higgs bundles.

For a parabolic $\U(p,q)$-Higgs bundle there is an invariant,
similar to the Toledo invariant in the non-parabolic case. We show
that this parabolic Toledo invariant has a  bound  provided by a
generalization of the Milnor--Wood inequality. Our main result in
the paper is to show that if the genus of the surface and the
number of marked points are both at least one, then the moduli
space of parabolic $\U(p,q)$-Higgs bundles with fixed topological
type, generic parabolic weights and full flags is non-empty and
connected if and only if the parabolic Toledo invariant satisfies
a generalized Milnor--Wood inequality (see Theorem
\ref{main-theorem}).

As in \cite{bgg} and \cite{ggm}, the main strategy  is to use the
Bott-Morse-theoretic techniques introduced by Hitchin \cite{h}.
The connectedness  properties of our  moduli space reduce to the
connectedness of a certain moduli space of parabolic triples
introduced in \cite{bg} in connection to the study of the
parabolic vortex equations and instantons of infinite energy. Much
of the paper is devoted to a thorough study of these moduli spaces
of triples and their connectedness properties.

After spelling out the correspondence between parabolic
$\U(p,q)$-Higgs bundles and representations of the fundamental
group of the punctured surface in $\U(p,q)$, we transfer our
results on connectedness of  the moduli space of parabolic
$\U(p,q)$-Higgs bundles to the moduli space of representations
(see Theorems \ref{nahodge-upq} and \ref{nahodge-upq2}). We then
apply this to the study of representations of the fundamental
group of certain complex elliptic surfaces of general type (see
Theorem \ref{elliptic-surface}). These are complex  surfaces whose
fundamental group is isomorphic to the orbifold fundamental group
of an orbifold Riemann surface.

We should point out that our main results do not apply when  the genus of the
Riemann surface is zero. This is not surprising if we have in mind that
on $\PP^1$ the parabolic weights must satisfy certain inequalites in order for
parabolic bundles to exist  (\cite{Bis,Bel}).
Presumably, something similar must be true also in the case of
parabolic $\U(p,q)$-Higgs bundles. We plan to come back to this
problem in a future paper.

In the process of finishing our paper we have come across several
papers  (\cite{bi,km,kr}) that seem to be related to our work in the
case of $\U(p,1)$. It would be interesting to investigate further
the relationship between these different approaches.

\vskip5pt

\noindent {\bf Acknowledments:} We thank the referee for a very
careful reading of the manuscript and for numerous suggestions.

\section{Parabolic Higgs bundles}\label{sec:PHB}

Let $X$ be a closed, connected, smooth  Riemann surface  of genus
$g\ge 0$ together with a finite set of marked points $x_{1},
\ldots, x_{s}$. Denote by $D$ the effective divisor
$D=x_{1}+\cdots +x_{s}$ defined by the marked points. A
\emph{parabolic vector bundle} $E$ over $X$ consists of a
holomorphic vector bundle together with a parabolic structure at
each $x\in D$, that is, a weighted flag on the fibre $E_{x}$,
  \begin{eqnarray*}
  & E_{x}=E_{x,1}\supset E_{x,2}\supset \cdots \supset E_{x,r(x)+1}=\{0\}, \\
  & 0\leq \alpha_{1}(x)< \ldots <\alpha_{r(x)}(x) <1.
  \end{eqnarray*}
We denote $k_i(x)=\dim (E_{x,i}/E_{x,i+1})$ the
\emph{multiplicity} of the weight $\alpha_{i}(x)$. It will
sometimes be convenient to repeat each weight according to its
multiplicity, i.e., we set $\tilde \alpha_1(x)=\ldots
=\tilde\alpha_{k_1(x)}(x)=\alpha_1(x)$, etc. We then have weights
$0\leq \tilde \alpha_{1}(x)\leq \ldots \leq\tilde \alpha_n(x) <1$,
where $n=\rk E$. Denote also
$\alpha(x)=(\tilde\alpha_{1}(x),\ldots , \tilde\alpha_{n}(x))$ the
system of weights at $x$ of $E$ and by $\alpha=(\alpha(x))_{x\in
D}$ the \emph{weight type} of $E$. We say that the flags are
\emph{full} if $k_{i}(x)=1$ for all $i$ and $x\in D$. Note that in
this case $\alpha(x)=(\tilde\alpha_{1}(x),\ldots ,
\tilde\alpha_{n}(x))=(\alpha_{1}(x),\ldots , \alpha_{n}(x))$. A
holomorphic map $f:E\to E'$ between parabolic bundles is called
\emph{parabolic} if $\alpha_{i}(x)>\alpha'_j(x)$ implies
$f(E_{x,i})\subset E'_{x,j+1}$  for all $x\in D$, and  $f$ is
\emph{strongly parabolic} if $\alpha_{i}(x)\ge \alpha'_j(x)$
implies $f(E_{x,i})\subset E'_{x,j+1}$ for all $x\in D$, where we
denote by $\alpha'_{j}(x)$ the weights on $E'$. We denote
$\PH(E,E')$ and $\SPH(E,E')$ the sheaves of parabolic and strongly
parabolic morphisms from $E$ to $E'$, respectively. If $E'=E$ we
denote these sheaves by $\PE(E)$ and $\SPE(E)$, respectively.

We define the \emph{parabolic degree} and \emph{parabolic slope}
of $E$ by
 \begin{eqnarray}\label{parabolic-degree}
 \pdeg(E)&=&\deg(E)+ \sum_{x\in
 D}\sum_{i=1}^{r(x)}k_i(x)\alpha_{i}(x), \\
 \pmu(E)&=&\frac{\pdeg(E)}{\rk(E)}.
 \end{eqnarray}
A parabolic bundle $E$ is said to be \emph{(semi)-stable} if for
every non-trivial proper parabolic subbundle $E'$ of $E$ we have
$\pmu(E')<\pmu(E)$ (resp.\ $\pmu(E')\le\pmu(E)$).

In the following we will use the following construction for
parabolic bundles, called \emph{parabolic direct sum}. Let $V$ and
$W$ two parabolic bundles with weight types $\alpha$ and $\alpha'$
we say that $E$ is the parabolic direct sum of $V$ and $W$ if and
only if $E=V\oplus W$ as holomorphic bundles, the system of
weights, $\tilde{\alpha}$, on $E$ consists of the ordered
collection of the weights in $\alpha$ and $\alpha'$, and the
corresponding filtration is such that
 $$
 E_{x,k}=V_{x,i}\oplus W_{x,j}
 $$
where $i$ (resp. $j$) is the smallest integer such that
$\tilde{\alpha}_{k}(x)\le \alpha_{i}(x)$ (resp.
$\tilde{\alpha}_{k}(x)\le \alpha'_{j}(x)$).

A \emph{parabolic Higgs bundle} is a pair $(E,\Phi)$ consisting of
a parabolic bundle $E$ and $\Phi\in H^{0}(\SPE(E)\otimes K(D))$,
i.e.\ $\Phi$ is a meromorphic endomorphism valued one-form with
simple poles along $D$ whose residue at $x\in D$ is nilpotent with
respect to the flag. A parabolic Higgs bundle is called
(semi)-stable if for every $\Phi$-invariant subbundle $E'$ of $E$,
its parabolic slope satisfies $\pmu(E')< \pmu(E)$ (resp.\
$\pmu(E)\leq \pmu(E)$), and  it is said to be  polystable if it is
the direct sum of stable parabolic Higgs bundles of the same
parabolic slope.

Fixing the topological invariants $n=\rk E$ and $d=\deg E$ and the
weight type $\alpha$, the moduli space $\cM=\cM(n,d;\alpha)$ is
defined as the set of isomorphism classes of polystable parabolic
Higgs bundles of type $(n,d;\alpha)$. Using Geometric Invariant
Theory, Yokogawa \cite{y1,y2} has showed that $\cM$ is a complex
quasi-projective variety, which is smooth at the stable points.

A \emph{parabolic $\U(p,q)$-Higgs bundle} on $X$ is a parabolic
Higgs bundle $(E,\Phi)$ such that $E=V\oplus W$, where $V$ and $W$
are parabolic vector bundles of ranks $p$ and $q$ respectively,
and
  $$
  \Phi=\left(\begin{array}{ll}  0 & \beta \\ \gamma & 0
  \end{array}\right) :(V\oplus W)\to (V\oplus W)\otimes K(D),
  $$
where  $\beta:W\to V\otimes K(D)$ and $\gamma: V\to W\otimes K(D)$
are strongly parabolic morphisms. A parabolic $\U(p,q)$-Higgs
bundle $(E=V\oplus W,\Phi)$ is (semi)-stable if the slope
stability condition $\pmu(E')<\pmu(E)$ (resp.\
$\pmu(E')\le\pmu(E)$) is satisfied for all $\Phi$-invariant
parabolic subbundles of the form $E'=V'\oplus W'$, i.e. for all
parabolic subbundles $V'\subset V$ and $W'\subset W$ such that
$\beta(W')\subseteq V'\otimes K(D)$ and $\gamma (V')\subseteq
W'\otimes K(D)$. Note that, a priori, this definition of stability
seems to be weaker than the stability definition for parabolic
Higgs bundles (we ask for $V'\subset V$ and $W'\subset W$). But
this is not the case, since for any $\Phi$-invariant $E'\subset
E$, we apply the $\U(p,q)$-stability condition to $V'\oplus W'$
and to $V''\oplus W''$, where $V'=V\cap E'$, $W'=W\cap E'$,
$V''=\pi_V(E')$, $W''=\pi_W(E')$ (where $\pi_V$, $\pi_W$ are the
projections of $V\oplus W$ onto $V$, $W$, respectively). Then
using the exact sequences $V'\to E'\to W''$ and $W'\to E'\to V''$,
one gets easily that $\pmu(E')\le \pmu(E)$).

Fix the topological invariants $a=\deg V$ and $b=\deg W$ and the
weight types $\alpha$ and $\alpha'$ for $V$ and $W$, respectively.
This determines a system of weights $\tilde \alpha$ and a flag structure, given by the parabolic direct sum construction, on
$E=V\oplus W$. Let
  $$
  \cU=\cU(p,q,a,b;\alpha,\alpha')
  $$
be the moduli space  of polystable  parabolic $\U(p,q)$-Higgs bundles of
degrees $(a,b)$ and weights  $(\alpha,\alpha')$.

We say that the weights are
\emph{generic} when every semistable parabolic Higgs bundle is
automatically stable, that is, there are no properly semistable
parabolic Higgs bundles. We will keep the following
assumption on the weights all throughout the paper (although some
of the results hold in more general situations):

\begin{assumption}\label{assumption}
The weights of $(E,\Phi)$ are
generic and  $(E,\Phi)$ has full flags at each parabolic point.
This means that all the weights of $V$ and $W$ are different and of
multiplicity one.
\end{assumption}

Note that the set of weights such that, for fixed degree and rank
of $E$, make $(E,\Phi)$ strictly semistable has positive
codimension. This justifies the term generic for the weights which
do not allow strict semistability.

The construction of $\cU$ follows the same arguments given in the
non-parabolic case (see \cite{bgg}).

\begin{prop} \label{prop:1}
Let $n=p+q$, $d=a+b$, and let  $\tilde\alpha$
be the system of weights defined by $\alpha$ and $\alpha'$ as above.
Then $\cU(p,q,a,b;\alpha,\alpha')$ embeds as a closed subvariety
in $\cM(n,d;\tilde\alpha)$.
\end{prop}

\begin{proof}
The proof is similar to that in the non parabolic case (see
Proposition 3.11 in \cite{bgg}). One only notices that in the case
$p=q$, the parabolic bundles $V$ and $W$ can not be parabolically
isomorphic since they have different weights.
\end{proof}

\begin{rmk} Sometimes we refer to elements $(E,\Phi)\in\cM$ as
parabolic $\GL(n,\CC)$-Higgs bundles, since  the structure group
of the frame bundle of $E$ is $\GL(n,\CC)$.
\end{rmk}

\section{Deformation theory} \label{sec:deformation}
The results of Yokogawa \cite{y1} and \cite{bgg}
readily adapt to describe the
deformation theory of parabolic $\U(p,q)$-Higgs bundles.

Let $(E = V \oplus W, \Phi)$ be a
parabolic $\U(p,q)$-Higgs bundle.  We introduce the following notation:
  \begin{align*}
    U & = \PE(E),              & \hat U & = \SPE(E),  \\
  U^+ &= \PE(V) \oplus \PE(W), &  \hat U^+ &= \SPE(V) \oplus \SPE(W), \\
U^- &= \PH(W,V)  \oplus \PH(V,W), &  \hat U^- &= \SPH(W,V)  \oplus \SPH(V,W).
\end{align*}

With this notation, $U=U^+\oplus U^-$, $\hat U=\hat U^+\oplus \hat U^-$,
$\Phi \in H^0(\hat U^-\otimes
K(D))$, and  $\ad(\Phi)$ sends  $U^+$ to $\hat U^-$ and
$U^-$ to $\hat U^+$. We consider the complex of sheaves
\begin{equation}
  \label{eq:tangentspace}
  C^{\bullet} :
U^+  \xrightarrow{\ad(\Phi)}  \hat U^- \otimes K(D).
\end{equation}

\begin{lem}
 \label{lem:stable-upq-vanishing}
  Let $(E,\Phi)$ be a stable parabolic $\U(p,q)$-Higgs bundle.  Then
\begin{align}
  \label{eq:ker-phi-u}
  \ker\bigl(\ad(\Phi) \colon H^0(U^+) \to
  H^0(\hat U^-\otimes K(D))\bigr) &= \CC, \\
  \label{eq:ker-phi-u-perp}
  \ker\bigl(\ad(\Phi) \colon H^0(U^-) \to
  H^0(\hat U^+\otimes K(D))\bigr) &= 0.
\end{align}
\end{lem}
\begin{proof}
Since $(E,\Phi)$ is stable as a
parabolic $\GL(n,\CC)$-Higgs bundle, it is simple, that is, its only
  endomorphisms are the non-zero scalars.  Thus,
$$
\ker\bigl(\ad(\Phi) \colon H^0(U) \to
H^0(\hat U\otimes K(D)) \bigr) = \CC.
$$
Since
$U = U^+ \oplus U^-$ and
$\ad(\Phi)$ sends $U^+$ to $\hat U^-$ and $U^-$ to $\hat U^+$,
the statements of the Lemma
follow.
\end{proof}

\begin{prop}
  \label{prop:upq-deformation}
\begin{itemize}
\item[{\rm (i)}] The space of endomorphisms of $(E,\Phi)$ is
isomorphic to
  the zeroth hypercohomology group $\HH^0(C^\bullet)$.
\item[{\rm (ii)}] The space of infinitesimal deformations of
$(E,\Phi)$ is
  isomorphic to the first hypercohomology group $\HH^1(C^\bullet)$.
\item[{\rm (iii)}] There is a long exact sequence
\begin{multline}  \label{eq:long-exact-tangent}
  0 \lto  \mathbb{H}^0(C^{\bullet}) \lto H^0(U^+) \lto
  H^0(\hat U^-\otimes K(D)) \lto  \mathbb{H}^1(C^{\bullet})  \\
  \lto  H^1(U^+) \lto  H^1(\hat U^-\otimes  K(D)) \lto
\mathbb{H}^2(C^{\bullet})\lto 0,
\end{multline}
where the maps $H^i(U^+) \lto  H^i(\hat U^-\otimes  K(D))$ are induced by
$\ad(\Phi)$.
\end{itemize}
\hfill\qed
\end{prop}

\begin{prop}
\label{prop:stable-upq-vanishing}
Let $(E,\Phi)$ be a stable parabolic  $\U(p,q)$-Higgs bundle, then
\begin{enumerate}
\item[{\rm (a)}] $\mathbb{H}^0(C^{\bullet}) = \CC$ (in other words
$(E,\Phi)$ is simple) and \item[{\rm (b)}]
$\mathbb{H}^2(C^{\bullet}) = 0$.
\end{enumerate}
\end{prop}
\begin{proof}
  (a) This follows immediately from
  Lemma~\ref{lem:stable-upq-vanishing} and (iii) of
  Proposition~\ref{prop:upq-deformation}.

  (b) For parabolic bundles $E$ and $F$
the sheaves $\PH(E,F))$ and $\SPH(F,E)\otimes \cO(D)$
are naturally dual to each other
(see for example \cite{by})
and  we thus have that
  $$
  \ad(\Phi) \colon H^1(U^+) \to H^1(\hat U^-\otimes K(D))
  $$
  is Serre dual to $\ad(\Phi) \colon H^0(U^-) \to  H^0(\hat U^+ \otimes K(D))$.
Hence Lemma~\ref{lem:stable-upq-vanishing} and (iii) of
  Proposition~\ref{prop:upq-deformation} show that
  $\mathbb{H}^2(C^{\bullet}) = 0$.
\end{proof}

\begin{prop} \label{prop:dim}
Assuming Assumption \ref{assumption}, the moduli space $\cU$ of stable parabolic $\U(p,q)$-Higgs
bundles is a smooth complex variety of dimension
 \begin{equation}\label{eqn:dimU}
 1+(g-1)(p+q)^{2}+\frac{s}{2}\big((p+q)^{2}-(p+q)\big),
 \end{equation}
where $g$ is the genus of $X$, and $s$ is the number of
marked points.
\end{prop}

\begin{rmk}
The formula in (\ref{eqn:dimU}) is also valid in the case $s=0$
and genus $g\ge 2$. In such case we recover the formula for the
dimension of the moduli space of non parabolic $\U(p,q)$-Higgs
bundles given in \cite{bgg}. As expected, this dimension is half
the dimension of the moduli space $\cM$ of parabolic
$\GL(n,\CC)$-Higgs bundles of rank $n=p+q$. Observe also that, in
order to have a non empty moduli space we need  $s\ge 3$ when
$g=0$.
\end{rmk}

\begin{proof}
Our assumption on the genericity of the weights implies that there
are no properly semistable parabolic $\U(p,q)$-Higgs bundles
and hence every point in $\cU$ is stable.
Smoothness follows from Propositions  \ref{prop:upq-deformation}
and \ref{prop:stable-upq-vanishing}. Now, our  assumption on having full
flags and different weights on $V$ and $W$ imply that
 $$
  \SPH(V,W)=\PH(V,W),
 $$
and
 \begin{eqnarray*}
 &&\dim\PH(V,W)_{x}+\dim\PH(W,V)_{x}= pq, \\
 &&\dim \PE(V)_x= \frac{p(p+1)}2,\\
 && \dim\PE(W)_{x}=\frac{q(q+1)}2.
 \end{eqnarray*}
Also, the short exact sequence
 $$
 0\to\PH(V,W)\to\Hom(V,W)\to \bigoplus_{x\in
 D}\frac{\Hom(V_{x},W_{x})}{\PH(V_{x},W_{x})}\to 0
 $$
implies that
 $$
 \deg(\PH(V,W))  = p\deg(W)-q\deg(V)+
  \sum_{x\in D}(\dim \PH(V_{x},W_{x})-pq).
 $$
Using the above information and Proposition \ref{prop:upq-deformation}
we have that the dimension
of the tangent space of $\cU$ at a point $(E,\Phi)$ is
 \begin{eqnarray*}
 \dim\HH^{1}(C^{\bullet}) & = &
 \dim\HH^{0}(C^{\bullet})
 + \dim\HH^{2}(C^{\bullet}) -\chi(C^{\bullet})\\
 &=&1-\chi(\PE(V)\oplus\PE(W))+\chi((\SPH(V,W)\oplus\SPH(W,V))\otimes K(D))\\
 &=&1-(p^{2}+q^{2})(1-g)-\deg(\PE(V))-\deg(\PE(W))
 +2pq((1-g)\\
 &&+\deg(\PH(V,W))+\deg(\PH(W,V))+ 2pq(2g-2)+2pqs\\
 &=&1+(g-1)(p+q)^{2}+2pqs+(p^{2}+q^{2}-2pq)s +\sum_{x\in
 D}\Big(\dim\PH(V,W)_{x} +\\
 &&+\dim\PH(W,V)_{x}-\dim\PE(V)_{x}-\dim\PE(W)_{x} \Big)\\
 &=&1+(g-1)(p+q)^{2}+\frac{s}{2}((p+q)^{2}-(p+q)).
 \end{eqnarray*}
\end{proof}

\section{Parabolic Toledo invariant} \label{sec:toledo}
In analogy with the non-parabolic case \cite{bgg}, one can associate
a Toledo invariant to a parabolic $\U(p,q)$-Higgs bundle.

\begin{defn}
The \emph{parabolic Toledo invariant} corresponding to the
parabolic Higgs bundle $(E=V\oplus W, \Phi)$ is
 \begin{equation}\label{eq:tau}
 \tau= 2\frac{pq}{p+q}(\pmu(V)-\pmu(W))
 \end{equation}
\end{defn}

The Toledo invariant will give us a way to classify components of
the moduli space of parabolic $\U(p,q)$-Higgs bundles. So we first
determine the possible values that it can take.

\begin{prop}\label{prop:ineq}
Let $(E=V\oplus W,\Phi=\left(\begin{array}{ll}  0 & \beta \\ \gamma & 0
  \end{array}\right))$ be a semistable
parabolic $\U(p,q)$-Higgs bundle. Then
 $$
 p(\pmu(V)-\pmu(E))\le \rk(\gamma)\left(g-1+\frac{s}{2}\right),
 $$
 $$
 q(\pmu(W)-\pmu(E))\le\rk(\beta)\left(g-1+\frac{s}{2}\right).
 $$
\end{prop}

\begin{proof}
Consider the parabolic bundles $N=\ker(\gamma)$ and
$I=\im(\gamma)\otimes K(D)^{-1}$. We have an exact
sequence of parabolic bundles
 $$
 0\to N\to V\to I\otimes K(D)\to 0
 $$
and
 \begin{equation}
 \label{eq:IKD}
 \begin{aligned}
 \pdeg(V)&=\pdeg(N)+\pdeg(I\otimes K(D)) \\
         &= \pdeg(N)+\pdeg(I)+\rk(I)(2g-2+s).
 \end{aligned}
 \end{equation}
Note that $I$ is a subsheaf of $W$ and the map $I\hookrightarrow
W$ is a parabolic map. Let $\tilde I\subset W$ be its saturation,
which is a subbundle of $W$, and endow it with the induced
parabolic structure. So $N$, $V\oplus \tilde I\subset E$ are
$\Phi$-invariant parabolic subbundles of $E$. The semistability of
$(E,\Phi)$ implies that
 \begin{equation} \label{eq:w1}
 \begin{aligned}
 \pmu(N)&\le&\pmu(E),\\
 \pmu(V\oplus I)\leq \pmu(V\oplus \tilde I)&\le&\pmu(E).
 \end{aligned}
 \end{equation}
This yields
 \begin{eqnarray*}
 \pdeg(N)&\le&\rk(N)\pmu(E),\\
 \pdeg(V)+\pdeg(I)&\le&(p+\rk(I))\pmu(E).
 \end{eqnarray*}
Adding both  and using (\ref{eq:IKD}) we have the
inequality
 $$
 2\pdeg(V)\le 2p \pmu(E)+\rk(I)(2g-2+s),
 $$
and hence
 $$
 p(\pmu(V)-\pmu(E))\le \rk(\gamma)\left(g-1+\frac{s}{2}\right).
 $$
 The other case is analogous.
\end{proof}

\begin{rmk} \label{rmk:sharp}
 The inequalities in Proposition \ref{prop:ineq} are not sharp.
 This is due to the fact that (\ref{eq:w1}) can be improved by
 assigning to $I$ the weights induced by the inclusion $I\subset W$.
\end{rmk}

One has the following bound for the Toledo invariant.

\begin{prop} \label{cor:bound_toledo}
Let $(E,\Phi)$ be a semistable parabolic $\U(p,q)$-Higgs
subbundle. Then,
 $$
 |\tau|\le \tau_M=\min\{p,q\}(2g-2+s),
 $$
\end{prop}

\begin{proof}
 Noting that
 \begin{equation}\label{eq:equs}
 \pmu(E)=\frac{p}{p+q}\pmu(V)+\frac{q}{p+q}\pmu(W),
 \end{equation}
 Proposition \ref{prop:ineq} may be rewritten as
 \begin{eqnarray*}
 q(\pmu(E)-\pmu(W))&\le&\rk(\gamma)\left(g-1+\frac{s}{2}\right),\\
 p(\pmu(E)-\pmu(V))&\le&\rk(\beta)\left(g-1+\frac{s}{2}\right).
 \end{eqnarray*}
 By (\ref{eq:equs}) we also have
 $\tau=2p(\pmu(V)-\pmu(E))=2q(\pmu(E)-\pmu(W))$. The result
 follows.
\end{proof}

\section{Hitchin equations and parabolic Higgs bundles}
\label{gauge}

In order to study the topology of $\mathcal{U}$ we  need a
gauge-theoretic interpretation of this moduli space in terms of
solutions to  the Hitchin  equations. One can adapt the arguments
given by  Simpson \cite{s2} for the case of parabolic
$\GL(n,\CC)$-Higgs bundles to the $\U(p,q)$ situation, along the
lines of what is done in \cite{bgg} in the non-parabolic case.
Similarly, to  construct  the moduli space from this point of
view, one can adapt the construction  given by Konno \cite{k} (see
also \cite{nst}) in the  parabolic $\GL(n,\CC)$ case.

A parabolic structure on a smooth vector bundle is defined in a
similar way to what is done in the holomorphic category. Let $E$
be a smooth parabolic vector bundle of rank $n$ and fix a
hermitian metric $h$ on  $E$ which is smooth in $X\setminus D$ and
whose (degenerate) behaviour around the marked points is given as
follows.  We say that a local frame $\{e_1,\dotsc,e_n\}$ for $E$
around $x$ \emph{respects the flag at $x$} if $E_{x,i}$ is spanned
by the vectors $\{e_{M_i+1}(x),\dotsc,e_n(x)\}$, where $M_i =
\sum_{j \leq i} k_j(x)$. Let $z$ be a local coordinate around $x$
such that $z(x) = 0$. We require that $h$ be of the form
\begin{displaymath}
  h =
  \begin{pmatrix}
    \abs{z}^{2\tilde{\alpha}_1} & & 0 \\
     & \ddots & \\
    0 & & \abs{z}^{2\tilde{\alpha}_n} \\
  \end{pmatrix}
\end{displaymath}
with respect to some local frame around $x$ which respects the
flag at $x$, where $\tilde{\alpha}_i=\tilde{\alpha}_i(x)$.

A  unitary connection $d_A$ associated to a
smooth $\dbar$ operator $\dbar_E$ on $E$  via the hermitian metric $h$
is singular at the marked points: if we write $z = \rho \exp(\sqrt{-1} \theta)$ and
$\{e_i\}$ is the local frame used in the definition of $h$, then
with respect to the local frame $\{\epsilon_i =
e_i/\abs{z}^{\tilde{\alpha}_i}\}$, the connection is of the form
 $$
  d_A = d + \sqrt{-1} \left(
  \begin{smallmatrix}
    \tilde{\alpha}_1 && 0 \\
    & \ddots & \\
    0 && \tilde{\alpha}_r
  \end{smallmatrix}
  \right)
  d\theta + A',
 $$
where  $A'$ is regular.
We denote the space of smooth $\dbar$-operators on
$E$ by ${\CCC}_E$, the space of associated $h$-unitary connections
by $\AAA_E$, the group of complex parabolic gauge
transformations by $\GGG^{\CC}_E$ and the subgroup of $h$-unitary
parabolic gauge transformations by $\GGG_E$.

Let $V$ and $W$ be smooth parabolic vector bundles  equipped with
hermitian metrics $h_V$ and $h_W$
adapted to the parabolic structures in the sense explained above.
We denote
$\CCC:=\CCC_V\times \CCC_W$,
$\GGG^\CC:=\GGG_V^\CC\times \GGG_W^\CC$,
$\GGG:=\GGG_V \times \GGG_W$. The
space of Higgs fields is
$\boldsymbol{\Omega} = \boldsymbol{\Omega}^+ \oplus  \boldsymbol{\Omega}^-$,
where
$\boldsymbol{\Omega}^+ =\Omega^{1,0}(\SPH(W,V)\otimes\cO(D))$ and
$\boldsymbol{\Omega}^- =\Omega^{1,0}(\SPH(V,W)\otimes\cO(D))$.
Here we  regard  $\SPH(W,V)$ and $\SPH(V,W)$ as smooth vector
bundles defined like in the holomorphic category.

Following Biquard \cite{b} and  Konno \cite{k}, we introduce
certain weighted Sobolev norms and  denote the corresponding
Sobolev completions of the spaces
defined above by ${\CCC}^{k}_{1}$, $\boldsymbol{\Omega}^k_1$,
${(\GGG^{\CC})}^k_2$ and ${\GGG}^k_2$.
Let
  $$
  {\HHH}= \{ (\dbar_E,\Phi) \in {\CCC} \times {\boldsymbol{\Omega}}
  \ | \  \dbar_E\Phi = 0 \}
  $$
and let ${\HHH}^k_1$ be the corresponding subspace of ${\CCC}^k_1
\times {\boldsymbol{\Omega}}^k_1$.

Let $\dbar_E=(\dbar_V,\dbar_W)$ where $\dbar_V\in \CCC_V$ and
$\dbar_W\in \CCC_W$, and $\Phi=\left(\begin{array}{ll}  0 & \beta \\ \gamma & 0
  \end{array}\right)$ with $\beta \in
\boldsymbol{\Omega}^+$ and $\gamma \in \boldsymbol{\Omega}^-$. Let
$F(A_V)$ and $F(A_W)$ be  the curvatures of the $h_V$ and
$h_W$-unitary connections corresponding to $\dbar_V$ and
$\dbar_W$, respectively. Let $\beta^*$ and $\gamma^*$ be the
adjoints with respect to $h_V$ and $h_W$. Fix a K\"ahler form $\omega$ on
$X$ with volume of $X$ normalized to $2\pi$. We consider the
moduli space $\cS$ defined by the subspace  of elements in
${\HHH}^k_1$ satisfying \emph{Hitchin equations}
\begin{align*}
    F(A_V) + \beta\beta^* + \gamma^*\gamma &=
    -\sqrt{-1}\mu\id_V\omega,  \\
    F(A_W) + \gamma\gamma^* + \beta^*\beta &=
    -\sqrt{-1}\mu\id_W\omega,
\end{align*}
modulo gauge transformations in  ${\GGG}^k_2$, where the equations
are only defined on $X\setminus D$. Taking the traces of the
equations, adding them,  integrating over $X\setminus D$, and
using the Chern--Weil formula for parabolic bundles, we find  that
$\mu=\pmu(V\oplus W)$.

The subspace of smooth points in ${\HHH}^k_1$ carries a K\"ahler metric
induced by  the complex structure of $X$ and the hermitian metrics $h_V$
and $h_W$. The Hitchin equations are moment map equations for the action
of ${\GGG}^k_2$ on  this subspace. In particular, the smooth part of
$\cS$, which corresponds to irreducible solutions, is obtained as a
K\"ahler quotient.
Under the genericity assumptions on the parabolic weights in
Assumption \ref{assumption}, all the solutions are irreducible and
the moduli space $\cS$ is a smooth K\"ahler manifold.

Fix the topological invariants $p=\rk V$, $q=\rk W$,
$a=\deg V$, $b=\deg W$ and the weight types $\alpha$ and $\alpha'$ of $V$
and $W$, respectively. Then
  $$
  \cU(p,q,a,b;\alpha,\alpha')\cong  {(\HHH^{s})}^k_1 / {(\GGG^{\CC})}^k_2,
  $$
where $\HHH^s$ are the stable elements in $\HHH$.
Moreover, if $\cS(p,q,a,b;\alpha,\alpha')$ is the moduli space of solutions
for these fixed invariants, we have the following.

\begin{thm}\label{hk}
There is a homeomorphism
  $$
  \cU(p,q,a,b;\alpha,\alpha')\cong  \cS(p,q,a,b;\alpha,\alpha').
$$
  \end{thm}

\section{Morse theory on the moduli space of parabolic $\U(p,q)$-Higgs
bundles.}\label{sec:morse}

In this section we recall the Bott-Morse theory used already in the study of
parabolic Higgs bundles in \cite{ggm,by}. There is an action of
$\CC^{\ast}$ on $\cU$ given by
 \begin{eqnarray*}
 \psi: \quad
 \CC^{\ast}\times \cU & \to& \cU \\
 (\lambda,(E,\Phi))&\mapsto& (E,\lambda\Phi).
 \end{eqnarray*}
This restricts  to a Hamiltonian action of the  circle  on the
moduli space $\cS$ of solutions to the Hitchin equations, which is
isomorphic to $\cU$ (Theorem \ref{hk}), with associated moment map
 $$
 [(E,\Phi)]\mapsto
 -\frac{1}{2}\|\Phi\|^{2}=-\imat \int_{X} \tr(\Phi\Phi^{\ast}).
 $$
We choose to  use the positive function, $f:\cU\to \RR$
 \begin{equation}\label{eq:f}
 f([E,\Phi])=\|\Phi\|^{2}.
 \end{equation}
Clearly $f$ is bounded below since it is non-negative. It is also
proper, this follows from the properness of the moment map
associated to the circle action on $\cM$ \cite{Bis} (see also
\cite{ggm}) and the fact that $\cU \subset \cM$ is a closed
subset.

To study the connectedness properties of $\cU$, we use the
following basic result: if $Z$ is a Hausdorff space and
$f:Z\to \RR$ is proper and bounded below then $f$ attains a
minimum on each connected component of $Z$. Therefore, if the
subspace of local minima of $f$ is connected then so is $Z$.
We thus have the following.

\begin{lem}\label{lem:conn}
The function $f:\cU\to \RR$ defined in {\rm (\ref{eq:f})} has a
minimum on each connected component of $\cU$. Moreover, if the
subspace of local minima of $f$ is connected then so is $\cU$.
\hfill $\Box$
\end{lem}

Now we will describe the minima of $f$. For this we introduce
the subset of $\cU$ defined by
\begin{equation}\label{minima}
 \cN=\cN(p,q,a,b;\alpha,\alpha') = \{(E,\Phi)\in \cU(p,q,a,b;\alpha,\alpha')
 \,\;\;\mbox{such that} \,\;\;  \beta=0\;\;\; \mbox{or}\;\;\; \gamma=0\}.
 \end{equation}

\begin{prop}
For every  $(E,\Phi)\in \cU$
 $$
 f(E,\Phi)\ge\frac{|\tau|}{2}\, ,
 $$
with equality if and only if $(E,\Phi)\in \cN$.
\end{prop}

\begin{proof}
The proof is similar to the one for Proposition 4.5 in \cite{bgg} apart from the fact
that we are using adapted metrics on the bundle.
\end{proof}

We will  prove that $\cN$ is the subvariety of local minima of $f$. For
this we have to describe the critical points of $f$ and characterize
the local minima. By a theorem of Frankel \cite{f}, the
critical points of $f$ are exactly the fixed points of the circle
action.


For a fixed point $(E,\Phi)$ of the circle action, we have an
isomorphism $(E,\Phi)\cong (E,e^{\imat\theta}\Phi)$ which yields
the following commutative diagram.
 \begin{displaymath}
 \begin{CD}
 E @>\Phi>> E\otimes K(D)\\
 @V\psi_{\theta}VV @VV\psi_{\theta}\otimes 1_{K(D)}V\\
 E @>e^{\imat\theta}\Phi>> E\otimes K(D).
 \end{CD}
 \end{displaymath}

\begin{prop}[{\cite[Thm.\ 8]{s2}}]\label{prop:simpson}
The equivalence class of a stable parabolic Higgs bundle
$(E,\Phi)$ is fixed under the action of $S^{1}$ if and only if it
is a \emph{parabolic Hodge bundle}. This
means that $E$ decomposes as  a direct sum
 $$
 E=E_{0}\oplus E_{1}\oplus \cdots \oplus E_{m}
 $$
of  parabolic bundles, such that $\Phi_{l}=\Phi|_{E_{l}}$ belongs
to $H^{0}(\SPH(E_{l},E_{l+1})\otimes K(D))$. If $\Phi_{l}\ne 0$,
then the weight of the isomorphism  $\psi_\theta: E
\longrightarrow E$ on $E_{l+1}$ is one plus the weight of
$\psi_\theta$ on $E_{l}$.
\end{prop}

The decomposition of $E$ is given by the eigenbundles
corresponding to the eigenvalues of the circle action on
$(E,\Phi)$.

\begin{cor}\label{cor:alt}
In the situation of Proposition \ref{prop:simpson}, if $(E,\Phi)$
is stable, then each $\Phi_{l}$ is nonzero and the $E_{l}$ are
alternately contained in $V$ and $W$.
\end{cor}

\begin{proof}
The proof goes similarly to the non parabolic case (see Proposition 4.10 from
\cite{bgg}).
\end{proof}

Now we want to compute the index of a critical point $(E,\Phi)$.
For this we need to write the complex in (\ref{eq:tangentspace}) in
terms of the eigenbundle decomposition provided by Proposition
\ref{prop:simpson}. Hence
 $$
 \PE(V)\oplus \PE(W)=\bigoplus_{-m\le 2k\le m}U_{2k}
 $$
 $$
 \SPH(V,W)\oplus\SPH(W,V)=\bigoplus_{-m\le 2k+1 \le m}\hat U_{2k+1}.
 $$
where
 \begin{equation}\label{eqn:u's}
 \begin{aligned}
 U_{l} &=\bigoplus_{i-j=l}\PH(E_{j},E_{i}),\\
 \hat U_{l} &=\bigoplus_{i-j=l}\SPH(E_{j},E_{i}).
 \end{aligned}
 \end{equation}
Therefore the deformation complex (\ref{eq:tangentspace})
for a parabolic $\U(p,q)$-Higgs bundle $(E,\Phi)$  can be written as
 $$
  C^\bullet \colon
  \bigoplus_{-m\le 2k \le m}U_{2k}\quad
  \stackrel{\ad(\Phi)}{\longrightarrow}
  \bigoplus_{-m\le 2k+1\le m}\hat U_{2k+1} \otimes K(D).
 $$
Each piece of this complex gives a subcomplex whose
hypercohomology gives an eigenspace of the tangent space
$T_{(E,\Phi)}\cU$ for the circle action.

\begin{prop}\label{prop:T}
Let $(E,\Phi)$ be a stable parabolic $\U(p,q)$-Higgs bundle which
represents a fixed point of the circle action on $\cU$. Then the
eigenspace of the Hessian of $f$ corresponding to the eigenvalue
$-2k$ is $\HH^{1}$ of the following complex
  $$
  C^{\bullet}_{2k}:U_{2k}\stackrel{\ad(\Phi)}{\longrightarrow}
  \hat U_{2k+1}\otimes K(D)\, .
  $$
\end{prop}

\begin{proof} Similar to the non parabolic case (see Proposition 4.11 from \cite{bgg}).
\end{proof}

\begin{cor}\label{cor:minimum}
 $(E,\Phi)$ is a local minimum of $f$ if and only if
 $\HH^{1}(C^{\bullet}_{2k})=0$ for all $k\ge 1$.
\end{cor}

\begin{prop}\label{prop:minimumcondition}
Let $(E,\Phi)$ be a stable parabolic $\U(p,q)$-Higgs bundle which
is a fixed point of the $S^{1}$-action on $\cU$. Then
$\chi(C^{\bullet}_{2k})\le 0$  for all $k\ge 1$, and
equality holds if and only if
 $$
 \ad(\Phi)|_{U_{2k}}:U_{2k}\to \hat U_{2k+1}\otimes K(D)
 $$
is an isomorphism of bundles.
\end{prop}

\begin{proof}
We want to get a bound for
 \begin{equation}\label{eqn:harto}
 \chi(C^{\bullet}_{2k})=\chi(U_{2k})-\chi(\hat U_{2k+1}\otimes
 K(D)).
 \end{equation}

The dual of each $U_{l}$ is
 \begin{eqnarray*}
 U_{l}^{\vee}&=&\bigoplus_{i-j=l}(\PH(E_{j},E_{i}))^{\vee}=
 \bigoplus_{i-j=l}\SPH(E_{i},E_{j}(D))=\hat U_{-l}(D).
 \end{eqnarray*}
The dual of $\ad(\Phi)|_{U_{2k}}$ is
 $$
 (\ad(\Phi)|_{U_{2k}})^t = \ad(\Phi)|_{U_{-2k-1}}\otimes 1_{K^{-1}}
 : U_{-2k-1} \otimes K^{-1} \to \hat U_{-2k}(D).
 $$

The vector bundle $\PE(E)$ has a natural parabolic structure
induced by the parabolic structure of $E$. In fact $\PE(E)$ as a
parabolic bundle is the parabolic tensor product of the parabolic
bundle $E$ and the parabolic dual of $E$ (see \cite{y1}), and
hence its parabolic degree is $0$. With respect to this parabolic
structure $(\PE(E),\ad(\Phi))$, where $\ad(\Phi): \PE(E)\to
\SPE(E)\otimes K(D)$, is a parabolic Higgs bundle. Now, the
stability of $(E,\Phi)$ implies the polystability of
$(\PE(E),\ad(\Phi))$. This can be seen by producing a solution to
the Hitchin equations on $(\PE(E),\ad(\Phi))$ out of the solution
on $(E,\Phi)$, which exists by Theorem \ref{hk}.  Since the
solution on $(\PE(E),\ad(\Phi))$ may not be irreducible, we only
have polystability (in particular, semistability) of
 $(\PE(E),\ad(\Phi))$.
The subbundles  $\ker(\ad(\Phi)|_{U_{2k}})$ and $\ker(\ad(\Phi)|_{U_{-2k-1}})$ of $\PE(E)$
are $\ad(\Phi)$-invariant  and hence we can apply the stability
condition on the parabolic slopes.
Since the ordinary  degree is smaller than the parabolic degree, we have
$\deg(\ker(\ad(\Phi)|_{U_{2k}}))\le 0$ and $\deg(\ker(\ad(\Phi)|_{U_{-2k-1}}))\le 0$.
Therefore we have the following chain of inequalities
\begin{eqnarray}\label{inequality-2k}
 \nonumber \deg(U_{2k}) & = & \deg(\ker(\ad(\Phi)|_{U_{2k}}))+ \deg(\im(\ad(\Phi)|_{U_{2k}})) \\ \nonumber
 & \le & \deg(\im(\ad(\Phi)|_{U_{2k}})) \\ \nonumber
 & \le & - \deg(\im((\ad(\Phi)|_{U_{2k}})^t)) \\ \nonumber
 &=& - \deg(\im(\ad(\Phi)|_{U_{-2k-1}} \otimes 1_{K^{-1}}))\\
  & =& - \deg(\im(\ad(\Phi)|_{U_{-2k-1}} )) + \rk(\im(\ad(\Phi)|_{U_{-2k-1}})) (2g-2)\\ \nonumber
 &= & \deg(\ker(\ad(\Phi)|_{U_{-2k-1}} )) - \deg (U_{-2k-1})+ \rk(\im(\ad(\Phi)|_{U_{-2k-1}}))
(2g-2)\\ \nonumber
 &\le & -\deg (U_{-2k-1})+ \rk(\im(\ad(\Phi)|_{U_{-2k-1}})) (2g-2) \\ \nonumber
 &=&\deg (\hat U_{2k+1}(D))+ \rk(\im(\ad(\Phi)|_{U_{-2k-1}})) (2g-2),
\end{eqnarray}
where we have used that $\rk(\im(h))=\rk(\im(h^{t}))$ and that
$\deg(\im(h))\le -\deg(\im (h^{t}))$ for any morphism of sheaves
$h$.

Using this we have that
\begin{eqnarray*}
 \chi(C^{\bullet}_{2k})&=&\deg(U_{2k})+\rk(U_{2k})(1-g)-\deg(U_{2k+1}\otimes
 K(D))-\rk(U_{2k+1})(1-g)\\
 &=&\deg(U_{2k})+\rk(U_{2k})(1-g)-\deg(U_{2k+1})-\rk(U_{2k+1})(g-1+s)\\
 &\le&\deg(\hat U_{2k+1}(D))+
 \rk(\im(\ad(\Phi)|_{U_{2k}}))(2g-2)+\rk(U_{2k})(1-g)-\deg(U_{2k+1})\\
 &&- \rk(U_{2k+1})(g-1+s) \\
 &=&(g-1)(2\rk(\im(\ad(\Phi)|_{U_{2k}}))-\rk(U_{2k})-\rk(U_{2k+1})),
\end{eqnarray*}
where we have used that $\hat U_{2k+1}=U_{2k+1}$ since all the weights
are different and of multiplicity $1$, and hence for $i\neq j$ it
is $\SPH(E_{i},E_{j})=\PH(E_{i},E_{j})$, since $E_i$ and $E_j$ are
different pieces in the decomposition of Proposition
\ref{prop:simpson}.
We thus have  $\chi(C^{\bullet}_{2k})\le 0$. If equality holds then
$\rk(\im(\ad(\Phi)|_{U_{2k}}))=\rk(U_{2k})=\rk(U_{2k+1})$, and also equality holds
in (\ref{inequality-2k}), showing that  $\ad(\Phi)|_{U_{2k}}$ is an isomorphism
as claimed.
\end{proof}

\begin{cor}\label{cor:localminimum}
Let $(E,\Phi)$ be a stable parabolic $\U(p,q)$-Higgs bundle which
represents a critical point of the Morse function $f$. This
critical point is a minimum if and only if
 $$
 \ad(\Phi)|_{U_{2k}}:U_{2k}\to \hat U_{2k+1}\otimes K(D)
 $$
is an isomorphism for all $k\ge 1$.
\end{cor}

\begin{proof}
By Corollary \ref{cor:minimum}, $(E,\Phi)$ is a local minimum if
and only if
 \begin{equation}\label{eq:condition1}
 \HH^{1}(C^{\bullet}_{2k})=0,\quad \forall k\ge 1.
 \end{equation}
Note that by Proposition \ref{prop:stable-upq-vanishing},
$\HH^0(C^{\bullet}_{2k})=0$ and
$\HH^2(C^{\bullet}_{2k})=0$, for $k\ge 1$. Hence
$(E,\Phi)$ is a local minimum if and only if
 $$
 \chi(C^{\bullet}_{2k}) = \sum (-1)^{i}\dim
 \HH^{i}(C^{\bullet}_{2k})= 0, \quad \forall k\ge 1.
 $$
By Proposition \ref{prop:minimumcondition}, this is equivalent to
requiring that
 $$
 \ad(\Phi):U_{2k}\to \hat U_{2k+1}\otimes K(D)
 $$
be an isomorphism of sheaves.
\end{proof}

Finally, we show that all these minima are in $\cN$.

\begin{prop} \label{prop:m2}
Let $(E,\Phi)=(E_{0}\oplus\cdots\oplus E_{m},\Phi)$ be stable and a fixed
point of the circle action, with $m\ge 2$. Then $(E,\Phi)$ is not
a local minimum.
\end{prop}

\begin{proof}
First note that $U_{l}=\hat{U}_{l}=0$ for $l>m$, and note also that for $l=m$, $U_{m}=\PH(E_{0},E_{m})$. Now we divide the proof conforming the different possibilities for $U_{l}$ and $\hat{U}_{l}$ as the number $m$ of terms in the bundle decomposition of $E$ is even or odd.

If $m$ is even then $2k=m$ and
 $$
 \ad(\Phi)|_{U_{m}}: \PH(E_{0},E_{m})\to 0
 $$
does not satisfy Corollary \ref{cor:localminimum}, hence
$(E,\Phi)$ is not a local minimum.

If $m\ge 2$ is odd, then $2k=m-1$ and
 $$
 \ad(\Phi)|_{U_{m-1}}:\PH(E_{0},E_{m-1})\oplus \PH(E_{1},E_{m})\to
 \SPH(E_{0},E_{m})\otimes K(D).
 $$
We will show that this is not an injective map of sheaves, and
therefore $(E_{0}\oplus\cdots \oplus E_{m},\Phi)$ is not a
minimum. We prove this in a small open set where all the bundles
trivialize. We need to find $\zeta=(\zeta_{1},\zeta_{2})\in
U_{m-1}$, $\zeta\ne 0$ such that $\ad(\Phi)|_{U_{m-1}}(\zeta)=0$, i.e.\ we
need to find $\zeta_{1}$ and $\zeta_{2}$ making the following
diagram commutative.
 \begin{displaymath}
 \begin{CD}
 E_{0}@>\Phi >> E_{1}\otimes K(D)\\
 @VV\zeta_{1}V @VV\zeta_{2}\otimes 1_{K(D)}V\\
 E_{m-1} @>\Phi >> E_{m}\otimes K(D)\\
 \end{CD}
 \end{displaymath}
For this, take $\zeta_{2}\ne 0$ such that $\zeta_{2}\otimes
1_{K(D)} (E_{1}\otimes K(D))\subset \Phi(E_{m-1})$, this is
possible by taking $\zeta_{2}$ as the composition of $\Phi_{l}$ in
Proposition \ref{prop:simpson} tensor the appropriate power of
$K(D)$, note that they are nonzero by Corollary \ref{cor:alt}. Now
take $\zeta_1$ such that
 $$
 \Phi\circ\zeta_{1}=(\zeta_{2}\otimes 1_{K(D)})\circ\Phi,
 $$
therefore $\Phi_{m-1}(\zeta)= (\zeta_{2}\otimes 1_{K(D)})
\circ\Phi -\Phi\circ\zeta_{1}=0$ with $\zeta\ne 0$. So
$\Phi_{m-1}$ is not injective.
\end{proof}

\begin{cor}\label{cor:conn2}
The subvariety of local  minima of $f:\cU(p,q,a,b;\alpha,\alpha')\to \RR$
coincides with the set  $\cN(p,q, a,b;\alpha,\alpha')$ defined in
(\ref{minima}).
\end{cor}

\begin{proof}
By Proposition \ref{prop:m2}, for $(E,\Phi)$ to be a minimum it
must have a decomposition of the form $E=E_{0}\oplus E_{1}$ with
$\Phi$ mapping $E_{0}$ into $E_{1}$. But by definition the only
possible decompositions are $E=V\oplus W$ with $\Phi=\left(\begin{array}{ll}  0 & 0 \\ \gamma & 0
  \end{array}\right)$ and
$E=W\oplus V$ with $\Phi=\left(\begin{array}{ll}  0 & \beta \\ 0 & 0
  \end{array}\right)$. So $(E,\Phi)\in \cN$.

Conversely, if $(E,\Phi)\in \cN$ then $m=1$ and
$U_{2k}=\hat U_{2k+1}=0$, for $k\ge 1$. So Corollary
\ref{cor:localminimum} applies and $(E,\Phi)$ is a minimum.
\end{proof}

Which of the two components of the Higgs field vanishes is given
by the following.

\begin{lem}\label{lem:gob}
Let $(E,\Phi)\in \cN$. Then the Toledo invariant $\tau\neq 0$ and
\begin{itemize}
 \item[(i)] $\gamma=0$ if and only if $\tau < 0$.
 \item[(ii)] $\beta=0$ if and only if $\tau > 0$.
\end{itemize}
\end{lem}

\begin{proof}
Observe that $\tau$ can not be equal to zero because this implies
$\gamma=\beta=0$ and then $(E,\Phi)$ cannot be stable.
The rest  follows directly  from the definition of the Toledo invariant.

\end{proof}

Our main goal in the rest of the paper is to show the following.

\begin{thm}\label{thm:connected-minima}
Suppose $g>0$. Then there is a value
 $$
 \tau_L = \min\{p,q\}(2g-2+s) - \frac{|p-q|}{p+q}\, \epsilon,
 $$
with $\epsilon>0$ explicitly computable (see Remark
\ref{rem:hartoxx}), such that the subvariety $\cN(p,q,
a,b;\alpha,\alpha')$ is non-empty and connected if and only if the
parabolic Toledo invariant $\tau$ satisfies
the bound 
$|\tau|\le \tau_L$. The moduli space of parabolic $\U(p,q)$-Higgs
bundles $\cU(p,q, a,b;\alpha,\alpha')$ is empty for
$|\tau|>\tau_L$.
\end{thm}

\begin{proof}
In the case $p\neq q$, the result will follow from Proposition
\ref{prop:iden} and Theorem \ref{thm:mainharto}. In the case
$p=q$, the result will follow from Propositions \ref{prop:iden}
and \ref{prop:harto11}, Corollary \ref{cor:irreducibility} and
Remark \ref{rem:irreducibility}. Note that $\tau_L=\tau_M$ for
$p=q$.
\end{proof}

Combining Theorem \ref{thm:connected-minima}, Corollary
\ref{cor:conn2} and Lemma \ref{lem:conn}, we have the main result
of our paper.

\begin{thm} \label{main-theorem}
Suppose $g>0$ and $s>0$. The moduli space of parabolic
$\U(p,q)$-Higgs bundles $\cU(p,q, a,b;\alpha,\alpha')$ is
non-empty and connected if and only if $|\tau|\le \tau_L$.The
moduli space is empty whenever $|\tau|>\tau_L$.
  \hfill $\Box$
\end{thm}

\begin{rmk}
 It is likely that Theorem \ref{main-theorem} holds more generally than
 under Assumption \ref{assumption}. It should be enough to assume
 that $V\oplus W$ have full flags, but arbitrary (non-generic) weights. The
 reason is that the assumption of full flags is strong enough to
 avoid the type of problem that comes up in Theorem 3.32 of
 \cite{bgg}, since all the weights are distinct.
 One way to prove this would be to show that the moduli spaces for
 different choices of weights
 are related by flips as with the moduli spaces of triples
 (as in \cite{th}).
\end{rmk}

\begin{rmk}
  Actually, in both Theorems \ref{thm:connected-minima} and
  \ref{main-theorem}, the case $|\tau|=\tau_L$ does not occur
  under Assumption \ref{assumption}. This is true since $\sigma=2g-2$ is
  not a critical value for the appropriate moduli space of triples
  appearing in Proposition \ref{prop:iden} (see Remark
  \ref{rem:prop-iden}). For $p=q$, it cannot happen that
  $|\tau|=\tau_M$, as pointed out in Remark
  \ref{rem:irreducibility}.
\end{rmk}

\section{Parabolic triples}\label{sec:ptriples}

In the previous section, we have concluded that it is necessary to
study the connectedness of the subspace $\cN$ of $\cU$. This
subset consists of  parabolic $\U(p,q)$-Higgs bundles with
$\gamma=0$ or $\beta=0$,  hence giving  rise in a natural way to
objects called  parabolic triples.

We recall the basics of parabolic triples from \cite{bg,ggm}. A
\emph{parabolic triple} is a holomorphic triple $T=(E_{1},
E_{2},\phi)$ where $E_{1}$ and $E_{2}$ are parabolic bundles and
$\phi:E_{2}\to E_{1}(D)$ is a strongly parabolic homomorphism,
i.e.\ $\phi\in H^{0}(\SPH(E_{2},E_{1}(D)))$. We denote
by $\alpha=(\alpha^{1},\alpha^{2})$  the parabolic system of weights for
the triple  $(E_1,E_2,\phi)$, where $\alpha^{i}$ is
the system of weights of $E_{i}$ with $i=1,2$.

For $\sigma\in \RR$ the parabolic \emph{$\sigma$-degree} and
\emph{$\sigma$-slope} of $T$ are defined as
 \begin{equation}\label{eq:sigmamu}
 \begin{array}{l}
 \pdeg_{\sigma}(T)=\pdeg(E_{1})+\pdeg(E_{2})+\sigma
 \rk(E_{2}),\\[5pt]
 \pmu_{\sigma}(T)={\displaystyle \frac{\pdeg E_{1}+\pdeg
 E_{2}}{\rk(E_{1})+\rk(E_{2})}+\sigma\frac{\rk(E_{2})}{\rk(E_{1})+\rk(E_{2})}}.
 \end{array}
 \end{equation}
A parabolic triple $T'=(E'_{1}, E'_{2},\phi')$ is a
\emph{parabolic subtriple} of $T=(E_{1}, E_{2},\phi)$ if
$E'_{i}\subset E_{i}$ are parabolic subbundles for $i=1,2$ and $\phi'=\phi|_{E'_{2}}$ being
$\phi(E'_{2})\subset E'_{1}(D)$. As usual, $T$ is called
$\sigma$-stable (resp.\ $\sigma$-semistable) if for any non-zero
proper subtriple $T'$ we have $\pmu_{\sigma}(T')<\pmu_{\sigma}(T)$
(resp.\ $\pmu_{\sigma}(T')\le \pmu_{\sigma}(T)$). The triple $T$
is called $\sigma$-polystable if it is the direct sum of parabolic
triples with the same parabolic $\sigma$-slope.

Let
 $$
 \cN_{\sigma}=\cN_{\sigma}(r_{1},r_{2},d_{1},d_{2};\alpha^1,\alpha^2)
 $$
be the moduli space of isomorphism classes of $\sigma$-polystable
triples with fixed system of weights $(\alpha^1,\alpha^2)$ and 
$r_{1}=\rk(E_{1})$, $r_{2}=\rk(E_{2})$, $d_{1}=\deg(E_{1})$,
$d_{2}=\deg(E_{2})$. Let
 $$
 \cN_{\sigma}^s\subset \cN_{\sigma}
 $$
be the open subset consisting of $\sigma$-stable triples.

\begin{prop}\label{prop:cotasigma}
A necessary condition for $\cN_{\sigma}(r_{1},r_{2},d_{1},d_{2},\alpha^1,
\alpha^2)$
to be non-empty is
\begin{eqnarray*}
 \sigma_{m}<\sigma<\sigma_{M}\quad &\mathrm{if}& r_{1}\ne r_{2}\\
 \sigma_{m}<\sigma\qquad \quad &\mathrm{if}& r_{1}=r_{2}
\end{eqnarray*}
where
\begin{eqnarray*}
 \sigma_{m}&=&\pmu(E_{1})-\pmu(E_{2})\\
 \sigma_{M}&=&\left(1+\frac{r_{1}+r_{2}}{|r_{1}-r_{2}|}\right)
 (\pmu(E_{1})-\pmu(E_{2}))+s
 \frac{r_{1}+r_{2}}{|r_{1}-r_{2}|},\qquad \text{if } r_{1}\ne r_{2}.
\end{eqnarray*}
\end{prop}

\begin{proof}
See Proposition 4.3 from \cite{ggm}.
\end{proof}

\begin{rmk}
We will see later on that there is an effective upper bound
$\sigma_L$ given by (\ref{eqn:sL}) which in general is  strictly
smaller than $\sigma_M$.
\end{rmk}

The correspondence between parabolic triples and
parabolic $\U(p,q)$-Higgs bundles goes as follows. Let $(E,\Phi)$
be a parabolic $\U(p,q)$-Higgs bundle with $\Phi=\beta: W\to
V\otimes K(D)$. This defines a triple $T=(E_{1},E_{2},\phi)$ where
$E_{1}=V\otimes K$, $E_{2}=W$, $\phi=\beta$. Conversely, given a
parabolic triple $T=(E_{1},E_{2},\phi)$ we get a
parabolic $\U(p,q)$-Higgs bundle with $\Phi=\left(\begin{array}{ll}  0 & \beta \\ 0 & 0
  \end{array}\right)$ by defining
$(E=V\oplus W, \Phi)$ where $V=E_{1}\otimes K^{-1}$, $W=E_{2}$ and
$\beta=\phi$. When $(E,\Phi)$ is a parabolic $\U(p,q)$-Higgs
bundle with $\Phi=\left(\begin{array}{ll}  0 & 0 \\ \gamma & 0
  \end{array}\right):V\to W\otimes K(D)$ we have an analogous
correspondence. That is, the corresponding triple to $(E,\Phi)$ is
$T=(W\otimes K, V,\gamma)$.

\begin{lem}\label{lem:corresp}
A parabolic $\U(p,q)$-Higgs bundle $(E,\Phi)$ with $\beta=0$ or
$\gamma =0$ is parabolically (semi)stable if and only if the
corresponding parabolic triple is $\sigma$-(semi)stable for
$\sigma=2g-2$.
\end{lem}

\begin{proof}
Let $T=(E_{1},E_{2},\phi)$ be the triple defined by $(E,\Phi)$
(without loss of generality we assume $\gamma=0$). Therefore if we
set $\sigma=2g-2$ we have
 \begin{eqnarray} \label{eq:pmut}
 \pmu_{\sigma}(T)&=&\frac{\pdeg(E_{1})+\pdeg(E_{2})}{\rk(E_{1})+\rk(E_{2})}+
 \sigma\frac{\rk(E_{2})}{\rk(E_{1})+\rk(E_{2})}\nonumber\\
 &=&\frac{\pdeg(V)+\pdeg(W)+p(2g-2)}{p+q}+ \sigma\frac{q}{p+q} \\
 &=& \pmu(E) + 2g-2 . \nonumber
 \end{eqnarray}
Note that the correspondence between parabolic triples and
$\U(p,q)$ parabolic bundles with $\beta=0$ or $\gamma=0$ gives
also a correspondence between parabolic subtriples and parabolic
subbundles. That is, given a subbtriple $T'$ of $T$ the
corresponding parabolic  $\U(p,q)$-Higgs bundle is a
$\Phi$-invariant subbundle of $(E,\Phi)$, and conversely given
$(E',\Phi')$ the corresponding triple gives a parabolic subtriple
of $T$. Hence equation (\ref{eq:pmut}) gives that
$\pmu_{2g-2}(T')<\pmu_{2g-2}(T)$ if and only if $\pmu(E')<\pmu(E)$
(analogously for the semistability condition).
\end{proof}

Combining  the arguments above and Lemma \ref{lem:gob}, we have
the following correspondence.

\begin{prop}\label{prop:iden}
Let $\cN(p,q,a,b;\alpha,\alpha')$ be the submanifold of local
minima of $\cU(p,q,a,b;\alpha,\alpha')$ and let $\tau$ be the
Toledo invariant then,
\begin{itemize}
 \item[(i)] If $\tau < 0$ then
$\cN(p,q,a,b;\alpha,\alpha')
=\cN_{2g-2}(p,q,a+p(2g-2),b;\alpha,\alpha')$.
 \item[(ii)] If $\tau > 0$ then
$\cN(p,q,a,b;\alpha,\alpha')
=\cN_{2g-2}(q,p,b+q(2g-2),a;\alpha',\alpha)$.
\end{itemize}
\end{prop}

\begin{proof}
It follows immediately from  Lemma \ref{lem:gob}.
\end{proof}

\begin{rmk} \label{rem:prop-iden}
 Note that the genericity condition on the weights implies that
 there are no properly $\sigma$-semistable triples for
 $\sigma=2g-2$, that is, $\cN_{2g-2}^s=\cN_{2g-2}$.
\end{rmk}

So we state the following assumption that we shall use during the
rest of the paper, and which is a translation of Assumption
\ref{assumption} via Proposition \ref{prop:iden}.

\begin{assumption}\label{assumption2}
We consider moduli spaces of $\sigma$-stable triples
$\cN_{\sigma}(r_{1},r_{2},d_{1},d_{2};\alpha^1,\alpha^2)$
satisfying that there are no properly $(2g-2)$-semistable triples
and such that all the weights are of multiplicity one, and the
weights of $E_1$ and $E_2$ are all different.
\end{assumption}

It is clear that in order for $\cN(p,q,a,b,\alpha,\alpha')$ to be
non-empty, $2g-2$ must be in the range for $\sigma$ given by
Proposition \ref{prop:cotasigma}, where  $\sigma_m$ and $\sigma_M$
are determined by the correspondence given in Proposition
\ref{prop:iden}. In fact, one has the following comparison of such
necessary condition with the Milnor--Wood inequality for the
parabolic Toledo invariant $\tau$ given in Proposition
\ref{cor:bound_toledo}

\begin{prop} \label{prop:harto11}
Let $\sigma_m$ and $\sigma_M$ be the bounds for $\sigma$ defined
in Proposition \ref{prop:cotasigma}  for the moduli space of
parabolic triples identified in Proposition \ref{prop:iden} with
the subvariety $\cN(p,q,a,b,\alpha,\alpha')$. Recall $\tau_M=\min\{p,q\}(2g-2+s)$. Then
 \begin{equation*}
 0\le|\tau|\le \tau_M\Leftrightarrow
 \begin{cases}
       \sigma_m\le 2g-2 \le
       \sigma_M, \qquad &\text{if $p\ne q$}, \\
       \sigma_m\le 2g-2, &\text{if $p=q$}.
 \end{cases}
 \end{equation*}
\end{prop}

\begin{proof}
 Write $\sigma_{m}$ and $\sigma_{M}$ in terms of $\tau$, that is,
 \begin{equation*}
 \begin{cases}
    \sigma_{m}=\frac{(p+q)}{2pq}\tau+2g-2, &\text{if $\tau<0$},\\
    \sigma_{m}=-\frac{(p+q)}{2pq}\tau+2g-2, &\text{if $\tau>0$},\\
    \sigma_{M}=\left(1+\frac{p+q}{|p-q|}\right)\left(\frac{(p+q)}{2pq}\tau+2g-2\right)
    +s\frac{p+q}{|p-q|}, &\text{if $\tau<0$},\\
    \sigma_{M}=\left(1+\frac{p+q}{|p-q|}\right)\left(-\frac{(p+q)}{2pq}\tau+2g-2\right)
    +s\frac{p+q}{|p-q|}, \qquad &\text{if $\tau>0$}.
 \end{cases}
\end{equation*}
{}From these equalities, the result is clear.
\end{proof}

\begin{rmk}
Proposition \ref{prop:harto11}
gives  a condition for the number of
marked points in order for  $\cN$ to be non-empty. Namely,
\begin{itemize}
 \item[(i)] If $g=0$ then $s\ge 3$,
 \item[(ii)] If $g=1$ then $s\ge 1$,
\end{itemize}
and no extra condition when  $g\ge 2$.
\end{rmk}

\section{Extensions and deformations of parabolic triples}

In order to study the differences between the moduli spaces
$\cN_{\sigma}$ as $\sigma$ changes, we need to study
extensions and deformations of parabolic triples. This study
is done in \cite{ggm}. We summarize the main results.

Let $T'=(E'_{1}, E'_{2},\phi')$ and $T''=(E''_{1},E''_{2},\phi'')$
be two parabolic triples. Let $\Hom(T'',T')$ denote the vector
space of homomorphisms from $T''$ to $T'$, and $\Ext^1(T'',T')$
be the vector space of extensions of the form
 $$
 0\to T'\to T\to T''\to 0,
 $$
that is, commutative diagrams:
\begin{displaymath}
 \begin{CD}
  0 @>>> E'_{2}@>>> E_{2} @>>> E''_{2} @>>> 0\\
  @.   @VV\phi'V @VV\phi V @VV\phi''V\\
  0@>>> E'_{1}(D)@>>> E_{1}(D) @>>> E''_{1}(D) @>>> 0.
 \end{CD}
\end{displaymath}

In order to study extensions of parabolic triples, we consider the
following complex of sheaves
 \begin{eqnarray}\label{eq:triplescompl}
 C^{\bullet}(T'',T'): \PH(E''_{1},E'_{1})\oplus \PH(E''_{2},
 E'_{2})&\to& \SPH(E''_{2},E'_{1}(D))\\
 (\psi_{1},\psi_{2})&\mapsto& \phi'\psi_{2}-\psi_{1}\phi''.\notag
 \end{eqnarray}

\begin{prop}[Proposition 4.7 \cite{ggm}]
There are natural isomorphisms
 \begin{eqnarray*}
 &&\Hom(T'',T')\cong \HH^{0}(C^{\bullet}(T'',T')),\\
 &&\Ext^{1}(T'',T')\cong \HH^{1}(C^{\bullet}(T'',T')).
 \end{eqnarray*}and a  long
exact sequence:
 \begin{equation}
 \begin{split}
 0&\to \HH^{0}\to H^{0}(\PH(E''_{1}, E'_{1})\oplus \PH(E''_{2},
 E'_{2}))\to H^{0}(\SPH(E''_{2},E'_{1}(D)))\\
 &\to \HH^{1}\to H^{1}(\PH(E''_{1}, E'_{1})\oplus \PH(E''_{2},
 E'_{2}))\to H^{1}(\SPH(E''_{2},E'_{1}(D)))\\
 &\to\HH^{2}\to 0.
 \end{split}
 \end{equation}
\end{prop}

We denote:
 \begin{equation}
 \begin{split}
 h^{i}(T'',T')=&\dim \HH^{i}(C^{\bullet}(T'', T')),\\
 \chi(T'',T')=& h^{0}(T'',T')-h^{1}(T'',T')+h^{2}(T'',T').
 \end{split}
 \end{equation}

\begin{prop}[Proposition 4.8 \cite{ggm}] For parabolic triples $T'$ and $T''$
 $$
 \chi(T'',T')=\chi(\PH(E''_{1},E'_{1}))+\chi(\PH(E''_{2},E'_{2}))
 -\chi(\SPH(E''_{2},E'_{1}(D))).
 $$
\end{prop}

\begin{cor}[Corollary 4.9 \cite{ggm}] \label{cor:exttriples}
For any extension $0\to T'\to T\to T''\to 0$ of parabolic triples
we have that
 $$
 \chi(T,T)=\chi(T',T')+\chi(T'',T'')+\chi(T'',T')+\chi(T',T'').
 $$
\end{cor}

Using the same arguments as in Proposition 3.5 of \cite{bgg2} one
can prove  the following.

\begin{prop}\label{prop:calch}
 Suppose that $T'$ and $T''$ are $\sigma$-semistable.
 \begin{itemize}
 \item[(i)]If $\pmu_{\sigma}(T')<\pmu_{\sigma}(T'')$,  then
 $\HH^{0}(C^{\bullet}(T'',T'))\cong 0$. \item[(ii)] If
 $\pmu_{\sigma}(T')=\pmu_{\sigma}(T'')$ and $T'$, $T''$ are $\sigma$-stable, then
  \begin{equation}\label{ec:calch}
  \HH^{0}(C^{\bullet}(T'',T'))\cong
   \begin{cases}
   \CC, \qquad \text{ if  $T'\cong T''$}\,,\\
   0 \,, \qquad \text{ if $T'\ncong T''$}.
   \end{cases}
  \end{equation}
 \end{itemize}
\end{prop}

\begin{thm}\label{th:triplesth}
Let $T=(E_{1}, E_{2},\phi)$ be a $\sigma$-stable parabolic triple.
\begin{itemize}
 \item[(i)] The Zariski tangent space at the point defined by $T$ in
  the moduli space $\cN_\sigma^s$ of $\sigma$-stable triples is isomorphic to
  $\HH^{1}(C^{\bullet}(T,T))$.
 \item[(ii)] If $\HH^{2}(C^{\bullet}(T,T))=0$, then the moduli
  space $\cN_\sigma^s$ of $\sigma$-stable parabolic triples is smooth in a
  neighbourhood of the point defined by $T$.
 \item[(iii)] $\HH^{2}(C^{\bullet}(T,T))=0$ if and only if the
  homomorphism
 $$
 H^{1}(\PE(E_{1}))\oplus H^{1}(\PE(E_{2}))\to
 H^{1}(\SPH(E_{2},E_{1}(D)))
 $$
  is surjective.
 \item[(iv)] At the smooth point in $\cN_{\sigma}^s$ represented by $T$,
  the dimension of the moduli space of $\sigma$-stable parabolic
  triples is
 \begin{equation*}
 \begin{split}\dim
 \cN_{\sigma}^s=&h^{1}(T,T)=1-\chi(T,T)\\
 =&1-\chi(\PE(E_{1}))-\chi(\PE(E_{2}))+\chi(\SPH(E_{2},E_{1}(D)))
 \end{split}
 \end{equation*}
 \item[(v)] If $\phi$ is injective or surjective then $T$ defines a
  smooth point in the moduli space.
\end{itemize}
\end{thm}
\begin{proof}
 The proof runs analogous to the non parabolic situation (see proof of Theorem 3.8 in \cite{bgg2}).
\end{proof}

\section{Critical values}\label{sec:critical}

A parabolic triple $T=(E_{1},E_{2},\phi)$ is \emph{strictly
$\sigma$-semistable} if and only if there is a proper subtriple
$T'= (E_{1}',E_{2}',\phi')$ such that
$\pmu_{\sigma}(T)=\pmu_{\sigma}(T')$, i.e.,
 \begin{equation}\label{eq:sigst}
  \pmu(T')+\sigma\frac{r_{2}'}{r_{1}'+r_{2}'}=\pmu(T)+\sigma
  \frac{r_{2}}{r_{1}+r_{2}},
 \end{equation}
where $r_1'=\rk (E_1')$, $r_2'=\rk(E_2')$. There are two ways in
which this can happen.  One is that there exists a parabolic
subtriple such that
 $$
 \frac{r_{2}'}{r_{1}'+r_{2}'}=\frac{r_{2}}{r_{1}+r_{2}}
 $$
therefore this implies
 $$
 \pmu(T')=\pmu(T).
 $$
In this case $T$ is strictly $\sigma$-semistable for all $\sigma$
(or at least for an interval of values of $\sigma$) and it is
called \emph{$\sigma$-independent semistable}. 
The other way in which strict $\sigma$-semistability can happen is
if equality holds for (\ref{eq:sigst}) but with
 $$
 \frac{r_{2}'}{r_{1}'+r_{2}'}\ne\frac{r_{2}}{r_{1}+r_{2}}.
 $$

\begin{defn} \label{def:critical}
The values of  $\sigma$ such that there exists a strictly
$\sigma$-semistable triple $T$ with a subtriple $T'$ such that
$\pmu_\sigma(T')=\pmu_\sigma(T)$ and
 $$
  \frac{r_{2}'}{r_{1}'+r_{2}'}\ne\frac{r_{2}}{r_{1}+r_{2}}
  $$
are called \emph{critical values}.
\end{defn}

\begin{prop}[Proposition 5.2 \cite{ggm}]\label{prop:modcrit}
\begin{itemize}
 \item[(i)] The critical values of $\sigma$ form a discrete subset
  of $[\sigma_{m},\sigma_{M}]$ if $r_1\neq r_2$, and of $[\sigma_m,\infty)$ if
  $r_1=r_2$.
 \item[(ii)] The stability criteria for two values of $\sigma$
  between two consecutive critical values are equivalent; thus the
  corresponding moduli spaces are isomorphic.
 \item[(iii)] For generic weights, $\sigma=2g-2$ is
  not a critical value.
\end{itemize}
\end{prop}

Let $\sigma_{c}$ be a critical value such that
$\sigma_{m}<\sigma_{c}<\sigma_{M}$. Here we adopt the convention
that  $\sigma_{M}=\infty$ when $r_{1}=r_{2}$. Set
 $$
 \sigma_{c}^{+}=\sigma_{c}+\epsilon, \quad
 \sigma_{c}^{-}=\sigma_{c}-\epsilon,
 $$
where $\epsilon>0$ is small enough so that $\sigma_{c}$ is the
only critical value in the interval
$(\sigma_{c}^{-},\sigma_{c}^{+})$.

\begin{lem}\label{lem:decmp}
Let $\sigma_{c}\in (\sigma_{m},\sigma_{M})$ be a critical value.
We define the flip loci $\cS_{\sigma_{c}^{\pm}}$ as the set of triples
in $\cN_{\sigma_{c}^{\pm}}^s$ which are $\sigma_{c}^{\pm}$-stable but not
$\sigma_{c}^{\mp}$-stable. Then
 $$
 \cN_{\sigma_{c}^{+}}^s-\cS_{\sigma_{c}^{+}}=
 \cN_{\sigma_{c}}^s=\cN_{\sigma_{c}^{-}}^s-\cS_{\sigma_{c}^{-}}.
 $$

\end{lem}

The following result is analogous to \cite[Proposition 5.4]{bgg2}.

\begin{prop}\label{prop:triplesext}
Let $\sigma_{c}\in(\sigma_{m},\sigma_{M})$ be a critical value.
Let $T=(E_{1},E_{2},\phi)$ be a triple which is
$\sigma_{c}$-semistable.
\begin{itemize}
\item[(1)] Suppose that $T$ represents a point in
$\mathcal{S}_{\sigma_c^+}$, i.e.\ suppose that $T$ is
$\sigma_c^+$-stable but not $\sigma_c^-$-stable.  Then $T$ has a
description as the middle term in an extension
 \begin{equation}\label{destab}
 0\to T'\to T\to T'' \to 0
 \end{equation}
in which
 \begin{itemize}
 \item[(a)]  $T'$ and $T''$ are both $\sigma_c^+$-stable, with
 $\pmu_{\sigma_c^+}(T')<\pmu_{\sigma_c^+}(T)$,
 \item[(b)] $T'$ and $T''$ are both $\sigma_c$-semistable with
 $\pmu_{\sigma_c}(T')=\pmu_{\sigma_c}(T)$.
 \end{itemize}
\item[(2)] Similarly, if $T$ represents a point in
$\mathcal{S}_{\sigma_c^-}$, i.e.\ if $T$ is $\sigma_c^-$-stable
but not $\sigma_c^+$-stable,  then $T$ has a description as the
middle term in an extension $\mathrm{(\ref{destab})}$ in which
 \begin{itemize}
 \item[(a)]  $T'$ and $T''$ are both $\sigma_c^-$-stable with
 $\pmu_{\sigma_c^-}(T')<\pmu_{\sigma_c^-}(T)$,
 \item[(b)] $T'$ and $T''$ are both $\sigma_c$-semistable with
 $\pmu_{\sigma_c}(T')=\pmu_{\sigma_c}(T)$.
 \end{itemize}
\end{itemize}
\end{prop}

The following lemma is proved with analogous arguments as in
Proposition 3.6 of \cite{bgg2}.

\begin{lem} \label{lem:W+}
Let $T'$ and $T''$ be triples which are $\sigma$-stable and of the
same $\sigma$-slope, for some $\sigma\geq 2g-2$. Then
 $$
 \HH^{2}(C^{\bullet}(T'',T'))=0.
 $$
\end{lem}

\begin{cor} \label{cor:expected}
 $\cN_\sigma$ is smooth of the expected dimension, for any
 $\sigma\ge 2g-2$.
\end{cor}

\begin{prop} \label{prop:codim}
 If $\sigma_{c} > 2g-2$ then the loci $\cS_{\sigma_{c}^{\pm}}\subset
 \cN_{\sigma_{c}^{\pm}}^s$ have codimension bigger than or equal to
 $-\chi(T',T'')$.
\end{prop}

\begin{proof}
Let us do the case of $\sigma_c^+$. For simplicity we denote
 \begin{eqnarray*}
 \cN'_{\sigma_{c}^{\pm}}&=&\cN_{\sigma_{c}^{\pm}}^s(r'_{1},
 r'_{2},d'_{1},d'_{2};\alpha^{1'},\alpha^{2'}),\\
 \cN''_{\sigma_{c}^{\pm}}&=&\cN_{\sigma_{c}^{\pm}}^s(r''_{1},
 r''_{2},d''_{1},d''_{2};\alpha^{1''},\alpha^{2''}).
 \end{eqnarray*}
It is known from \cite{y2} that $\cN'_{\sigma_{c}^\pm}$ and
$\cN''_{\sigma_{c}^\pm}$ are fine moduli spaces. That is, there are
universal parabolic triples $\cT'=(\cE'_{1},\cE'_{2},\Phi')$ and
$\cT''=(\cE''_{1},\cE''_{2},\Phi)$ over $\cN'_{\sigma_{c}^+}\times
X$ and $\cN''_{\sigma_{c}^+}\times X$ respectively. Thus we
consider the complex $C^{\bullet}(\cT'',\cT')$ as defined in
(\ref{eq:triplescompl}) and take relative hypercohomology with
respect to the projection
 $$
 \pi:X\times \cN'_{\sigma_{c}^+}\times \cN''_{\sigma_{c}^+}
 \to\cN'_{\sigma_{c}^+}\times \cN''_{\sigma_{c}^+}.
 $$
We define $W^{+}:=\HH^{1}_{\pi}(C^{\bullet}(\cT'',\cT'))$. By
Proposition \ref{prop:triplesext}, $\cS_{\sigma_{c}^+}$ is a
subset of the projective fibration $\PP W^+$ over
$\cN'_{\sigma_{c}^+}\times \cN''_{\sigma_{c}^+}$. The fibres of
this fibration are projective spaces of dimension
 $$
 \begin{aligned}
 \dim \PP(\Ext^{1}(T'',T'))=&\dim \Ext^{1}(T'',T')-1\\
 =&h^{0}(T'',T')+h^{2}(T'',T')-\chi(T'',T')-1\\
 =&-\chi(T'',T')-1,
 \end{aligned}
 $$
using Lemma \ref{lem:W+} and Proposition \ref{prop:calch} to
substitute $h^{0}(T'',T')=h^{2}(T'',T')=0$. Therefore
 $$
 \begin{aligned}
 \dim \cS_{\sigma_{c}^+} \leq & -\chi(T'',T')+
 \dim(\cN'_{\sigma_{c}^+}\times\cN''_{\sigma_{c}^+})\\
 =&-\chi(T'',T')-1+1-\chi(T',T')+1-\chi(T'',T'')\\
 =&\dim\cN_{\sigma_{c}^+}+\chi(T',T''),
 \end{aligned}
 $$
since the moduli spaces $\cN'_{\sigma_{c}^+}$ and
$\cN''_{\sigma_{c}^+}$ are smooth of the expected dimension.
Therefore $\dim\cN_{\sigma_{c}^+}^s-\dim\cS_{\sigma_{c}^+} \geq
-\chi(T',T'')$.
\end{proof}

Hence, if we prove that this codimension is positive then the
moduli spaces $\cN_{\sigma}^s$ for different values of $\sigma\geq
2g-2$ are birational, and in particular have the same number of
irreducible components.

\section{Codimension of the flip loci}\label{sec:chi}

Let  $\sigma_c$ be  a critical value in the interval
$(\sigma_m,\sigma_M)$ such that $\sigma_c\geq 2g-2$. Let $T'$ and
$T''$ be two $\sigma_c^\pm$-stable (and $\sigma_c$-semistable)
parabolic triples with $\pmu_{\sigma_c}(T')=\pmu_{\sigma_c}(T'')$.
Changing the roles of $T'$ and $T''$, we may compute the bound
$\chi(T'', T')$ for the codimension of the flip locus (Proposition
\ref{prop:codim}) using the complex (\ref{eq:triplescompl}). Under
our Assumption \ref{assumption2}, we have
$\SPH(E_2'',E_1'(D))=\PH(E_2'',E_1'(D))$, and hence the complex
(\ref{eq:triplescompl}) is
 $$
  \begin{array}{ccc}
  C^{\bullet}(T'',T'):C_{1}=\PH(E''_{1},E'_{1})\oplus
  \PH(E''_{2},E'_{2})&\stackrel{a_1}{\longrightarrow}
  & C_{0}(D)=\PH(E''_{2},E'_{1}(D)) \\
  (\xi_1,\xi_2) &\mapsto & \phi'\xi_2-\xi_1\phi''.
  \end{array}
 $$
Our task is  to bound the Euler characteristic of the complex
$C^\bullet(T'',T')$, that is,
 $$
 \chi(C^{\bullet}(T'',T'))=(1-g)(\rk(C_{1})-\rk(C_{0}))+
 \deg(C_{1})-\deg(C_{0}(D)).
 $$

In order to obtain  bounds for $\deg(C_{1})$ and $\deg(C_{0})$,
we follow a similar strategy to that used in \cite{bgg} in the non-parabolic
case, exploiting  the existence theorem for parabolic vortex
equations.

\begin{thm}[{\cite[Thm.\ 3.4]{bg}}] \label{thm:bg}
 Let $T=(E_{1},E_{2},\phi)$ be a
 parabolic triple. Let $\tau_{1}$ and $\tau_{2}$ satisfy
 $\tau_{1}\rk(E_{1})+\tau_{2}\rk(E_{2})=\pdeg(E_{1})+\pdeg(E_{2})$,
 and let $\sigma=\tau_{1}-\tau_{2}$. Then $E_{1}$ and $E_{2}$ admit
 hermitian metrics, adapted to the parabolic structures, satisfying
 \begin{eqnarray*}
 &&\imat \Lambda F(E_{1})+\phi\phi^{\ast}=\tau_{1}\id_{E_{1}},\\
 &&\imat \Lambda F(E_{2})-\phi^{\ast}\phi=\tau_{2}\id_{E_{2}},
 \end{eqnarray*}
 if and only if $T$ is $\sigma$-polystable. Here $F(E_i)$ is  the curvature
 of the hermitian metric of $E_i$ and $\Lambda$ is the contraction
with a K\"ahler form on $X$ with volume normalized to $2\pi$.
\end{thm}

One can easily show that
 \begin{eqnarray*}
 \tau_1 &=& \pmu_\sigma(T), \\
 \tau_2 &=& \pmu_\sigma(T)-\sigma.
 \end{eqnarray*}
Moreover, adding up the equations in Theorem \ref{thm:bg},
integrating, and using the Chern--Weil formula for parabolic bundles,
we have that
 $$
 r_{1}\tau_{1}+r_{2}\tau_{2}= \pdeg(E_{1}) + \pdeg (E_2).
 $$

In our situation,  the triples $T'$ and $T''$ are $\sigma$-stable
for $\sigma=\sigma_c^\pm$, and hence, by Theorem \ref{thm:bg},
there exist adapted hermitian metrics  such that
 \begin{eqnarray*}
 \imat \Lambda
 F(E_{1}')+\phi'(\phi')^{\ast}=\tau_{1}'\id_{E_{1}'},
 \,\,\,\,
 &&\imat \Lambda F(E_{2}')-(\phi')^{\ast}\phi'=\tau_{2}'\id_{E_{2}'}, \\
 \imat \Lambda
 F(E_{1}'')+\phi''(\phi'')^{\ast}=\tau_{1}''\id_{E_{1}''},
 &&\imat \Lambda F(E_{2}'')-(\phi'')^{\ast}\phi''=\tau_{2}''\id_{E_{2}''},
 \end{eqnarray*}
where $\sigma=\tau'_{1}-\tau'_{2}=\tau''_{1}-\tau''_{2}$. In
particular, $\tau'_{1}-\tau''_{1}=\tau'_{2}-\tau''_{2}$.

Let us consider the induced adapted hermitian metrics on $C_{0}$ and
$C_{1}$. The corresponding  curvatures are given by
 \begin{eqnarray*}
  F(C_{0})&=&-F(E''_{2})^{t}\otimes \id_{E'_{1}}+\id\otimes F(E'_{1}),\\
  F(C_{1})&=&\left(-F(E''_{1})^{t}\otimes \id_{E'_{1}}+\id_{E''_{1}}\otimes
  F(E'_{1}),-F(E''_{2})^{t}\otimes\id_{E'_{2}}+\id_{E''_{2}}\otimes
  F(E'_{2})\right).
 \end{eqnarray*}
Actually, we have defined $C_0$ and $C_1$ as holomorphic bundles,
but they admit parabolic structures in a natural way: given
parabolic bundles $E$ and $F$, there are parabolic duals $E^{*p}$
and parabolic tensor products $E\otimes^p F$ (see \cite{y1}
\cite{ggm}). Then the parabolic structure
on $\PH(E,F)$ is given by $E^{*p}\otimes^p F$. In the
formulas for $F(C_0)$ and $F(C_1)$ we have to consider the adapted
metrics for the parabolic structures on each $(E''_{j})^{*p}
\otimes^{p} E'_{i}$, induced by the adapted metrics on the bundles
$E'_{k}$ and $E''_{k}$, for $k=1,2$.

Consider the  homomorphism $a_{2}$ defined by
 \begin{eqnarray*}
 \PH(E''_{1},E'_{2})(-D)&\stackrel{a_{2}}{\longrightarrow}
 &\PH(E''_{1},E'_{1})\oplus \PH(E''_{2},E'_{2})\\
 \xi &\to&(\phi' \xi,\xi\phi'').
 \end{eqnarray*}
The connections on $C_{0}$ and $C_{1}$ satisfy
 \begin{equation}\label{eq:cvort}
 \begin{array}{l}
 \imat \Lambda
 F(C_{0})+a_{1}a_{1}^{\ast}=(\tau'_{1}-\tau''_{2})\id_{C_{0}}\\
 \imat \Lambda
 F(C_{1})-a_{1}^{\ast}a_{1}+a_{2}a_{2}^{\ast}=
 (\tau'_{1}-\tau''_{1})\id_{C_{1}}.
 \end{array}
 \end{equation}

\begin{lem}\label{lem:KQ}
Let $K$ and $Q(D)$ denote the kernel and the torsion-free part of
the cokernel, respectively, of the homomorphism $a_{1}$. Then
 \begin{eqnarray*}
 \pmu(K)&\le& \pmu_{\sigma}(T')-\pmu_{\sigma}(T''), \\
 \pmu(Q)&\ge& \pmu_{\sigma}(T'')-\pmu_{\sigma}(T')+\sigma.
 \end{eqnarray*}
\end{lem}

\begin{proof}
 The kernel $K$ is a subbundle of the hermitian bundle
 $C_{1}$, so that we may take the $\cC^{\infty}$ orthogonal
 splitting $C_{1}=K\oplus S$. Since $K$ is a holomorphic subbundle,
 the induced connection $D_K$ on $K$ satisfies
 $D_{C_{1}}|_{K}=D_{K}+A$, where $D_{C_{1}}$ is the connection on
 $C_{1}$ and $A\in \Omega^{1,0}(\Hom(K,S))$ is the second
 fundamental form of $K\subset C_{1}$. Therefore the curvature $F(K)$
 of the connection on $K$ satisfies
 $F(C_{1})|_{K}=F(K)+\bar{A}^{t}\wedge A$.

We now use the second equation in (\ref{eq:cvort}) restricted to
$K$, take the trace and integrate on $X\setminus D$, to get
 $$
 \int_{X\setminus D}\tr(\imat \Lambda (F(K)+\bar{A}^{t}\wedge
 A)-a_{1}^{\ast}a_{1}|_{K}+a_{2}a_{2}^{\ast}|_{K})=
 \int_{X\setminus D}\tr((\tau'_{1}-\tau''_{1})\id_{C_{1}}|_{K}).
 $$
That is
 $$
 \pdeg(K)+\|A\|^{2}_{L^{2}}+\int_{X\setminus D}\tr
 (a_{2}a_{2}^{\ast}|_{K})=(\tau'_{1}-\tau''_{1})\rk(K),
 $$
obtaining
 $$
 \pdeg(K)\le (\tau'_{1}-\tau''_{1})\rk(K)
 $$
as desired, since $\tau_1'=\pmu_\sigma(T')$ and
$\tau_1''=\pmu_\sigma(T'')$.

To get the second inequality, let $S'(D)$  be the saturation of
the image of $a_1$, which is holomorphic subbundle of $C_0(D)$.
Then there is a $\cC^{\infty}$ orthogonal splitting
$C_{0}=S'\oplus Q$. The curvature of the induced connection
on $Q$ satisfies $F(C_{0})|_{Q}=F(Q)+ B\wedge\bar{B}^{t}$ with
$B\in\Omega^{0,1}(\Hom(Q,S'))$. If we consider the first equation
in (\ref{eq:cvort}) restricted to $Q$, take the  trace and integrate,
we get
 \begin{equation*}
 \int_{X\setminus D}\tr(\imat \Lambda
 (F(Q)+B\wedge\bar{B}^{t})+a_{1}a_{1}^{\ast})|_{Q}=
 \int_{X\setminus D}\tr((\tau'_{1}-\tau''_{2})\id_{C_{0}}|_{Q}).
 \end{equation*}
That is,
 $$
 \pdeg(Q)-\|B\|^{2}_{L^{2}}=(\tau'_{1}-\tau''_{2})
 (\rk(C_{0})-\rk(\im(a_{1})).
 $$
Hence,
 \begin{equation}\label{eq:Q}
 \pdeg(Q)\ge (\tau'_{1}-\tau''_{2})(\rk(C_{0})-\rk(\im(a_{1})),
 \end{equation}
as stated.
\end{proof}

\begin{thm}\label{thm:codim}
 Let $T'$ and $T''$ be $\sigma_c^\pm$-stable parabolic triples over a punctured
 Riemann surface of genus $g>0$ such that
 $\pmu_{\sigma_c}(T')=\pmu_{\sigma_c}(T'')$ for $\sigma_c\geq 2g-2$.
 Suppose that the morphism $a_1$ is not an isomorphism of bundles.
 Then
  $$
  \chi(C^{\bullet}(T'',T'))<0.
  $$
\end{thm}

\begin{proof}
We have
 \begin{equation}\label{eqn:harto2}
 \begin{array}{rcl}
 \chi(C^{\bullet}(T'',T'))&=&(1-g)(\rk(C_{1})-\rk(C_{0}))+\deg(C_{1})-\deg(C_{0}(D))\\
 &=&(1-g)(\rk(C_{1})-\rk(C_{0}))+\deg(K)+\deg(\im(a_{1}))-\deg(C_{0}(D))\\
 &\leq &(1-g)(\rk(C_{1})-\rk(C_{0}))+\deg(K)-\deg(Q)\\
 &=&(1-g)(\rk(C_{1})-\rk(C_{0}))+\deg(K)-\deg(Q(-D)(D)).
 \end{array}
 \end{equation}
Observe that for any (non-zero) parabolic bundle $E$, $\deg(E(D))>
\pdeg(E)\ge \deg(E)$, where the strict inequality is given by the
fact that the weights on $E$ always satisfy $0\leq \alpha_{i}(x)
<1$, for all $i$ and all $x\in D$. Using this, the hypothesis
$\sigma\geq 2g-2$, and Lemma \ref{lem:KQ}, we have
 \begin{equation}\label{eqn:harto3}
 \begin{array}{rcl}
  \chi(C^{\bullet}(T'',T'))& \leq &(1-g)(\rk(C_{1})-\rk(C_{0}))+\pdeg(K)-\pdeg(Q(-D))\\
 &=&(1-g)(\rk(C_{1})-\rk(C_{0}))-\sigma(\rk(C_{0}(D))-\rk(\im(a_{1}))\\
 &\le&(1-g)(\rk(C_{1})-\rk(C_{0}))+2(1-g)(\rk(C_{0})-\rk(\im(a_{1}))\\
 &=&(1-g)(\rk(C_{1})+\rk(C_{0})-2\rk(\im(a_{1})) \\ &\leq & 0,
 \end{array}
 \end{equation}
using that $g\geq 1$. If either $K$ or $Q$ is a non-zero bundle,
then the first line of (\ref{eqn:harto3}) is a strict inequality.
If both are zero and $a_1$ is not an isomorphism, then the third
line of (\ref{eqn:harto2}) is a strict inequality since
$\im(a_1)\neq C_0(D)$. In both cases,
 $$
 \chi(C^{\bullet}(T'',T'))< 0.
 $$
\end{proof}

\begin{rmk} \label{rem:g=0}
Note that this theorem does not cover the case $g=0$.
This is not so surprising if we recall that, in order for
parabolic bundles to exist on $\PP^1$, the parabolic weights
must satisfy certain inequalites (\cite{Bis,Bel}).
Presumably, something similar must be true also  in the case of
parabolic $\U(p,q)$-Higgs bundles.
\end{rmk}

The following result will be useful in the next sections.

\begin{lem} \label{lem:geniso}
 If $a_1$ is generically an isomorphism of bundles, then either
 \begin{itemize}
  \item[(a)] $E''_1=0$ and $\phi':E_2'\to E_1'$ is generically an
  isomorphism. In this case, $r_2 >r_1$.
  \item[(b)] $E'_2=0$ and $\phi'':E_2''\to E_1''$ is generically an
  isomorphism. In this case, $r_2 <r_1$.
 \end{itemize}
\end{lem}

\begin{proof}
  One may look at a generic point $x\in X\setminus D$, i.e, a point
  where the maps $\phi'$ and $\phi''$ are generic. We have
 $$
  \begin{array}{ccc}
  (a_1)_x: \PH(\CC^{r_1''},\CC^{r'_1})\oplus \PH(\CC^{r_2''},\CC^{r'_2}) &\to
  &\PH(\CC^{r_2''},\CC^{r'_1})\\
  (\alpha,\beta) &\mapsto & \phi'_x \ \beta- \alpha \ \phi''_x\, .
  \end{array}
 $$

If $\phi''_x$ is not surjective, take $\beta=0$ and $\alpha\neq 0$
with $\alpha|_{\im(\phi''_x)}=0$. Then $(a_1)_x(\alpha,\beta)=0$.
If $\phi'_x$ is not injective, take $\alpha=0$ and $\beta\neq 0$
with $\im(\beta)\subset \ker\phi'_x$, to get
$(a_1)_x(\alpha,\beta)=0$. Both possibilities contradict the
injectivity of $(a_1)_x$. Therefore $\phi''_x$ is surjective and
$\phi'_x$ is injective.

If neither of $\phi'_x$ and $\phi''_x$ is an isomorphism, then
take a map $\CC^{r_2''}\to \CC^{r'_1}$ which induces a non-zero
map $\ker(\phi''_x) \to \coker (\phi'_x)$. This cannot be in the
image of $(a_1)_x$, contradicting our assumption. So either
$\phi'_x$ or $\phi''_x$ are isomorphisms. In the first case
$r'_1r''_1+r'_2r''_2=r''_2r'_1$ gives $r''_1=0$ and we are in case
(a). In the second, we are in case (b).
\end{proof}

\section{Irreducibility of the moduli space of triples for $r_1\neq r_2$}
\label{sec:r-different}

This section is devoted to study the irreducibility and
non-emptiness of the moduli space of $\sigma$-stable parabolic
triples for ranks $r_1\neq r_2$.

Given a triple $T=(E_{1},E_{2},\phi)$ one has the dual triple
$T^{\ast}=(E_{2}^{\ast p},E_{1}^{\ast p},\phi^{t})$ where
$E_{i}^{\ast p}$ is the parabolic dual of $E_{i}$ and $\phi^{t}$
is the dual of $\phi$.

\begin{prop}  \label{prop:duality}
The $\sigma$-stability of $T$ is equivalent to the
$\sigma$-stability of $T^{\ast}$. The map $T\mapsto T^{\ast}$
defines an isomorphism of the corresponding moduli spaces of
$\sigma$-stable triples.
\end{prop}

This allows us to restrict to the case $r_{1}> r_{2}$ and appeal
to duality for the case $r_{1}<r_{2}$. So throughout this section
we assume that $r_{1}>r_{2}$.

\begin{lem} \label{lem:extra}
  Let $X$ be a Riemann surface with a finite number
of marked points  and let $E$, $F$ be
  parabolic bundles on $X$. Let $p\in X$ be a parabolic point.
 Then there is a natural exact sequence
  $$
   0 \to \frac{\Hom(E_p,F_p)}{\PH(E_p,F_p)} \otimes \cO(-p)
   \to \PH(E,F)_p \to \PH(E_p,F_p) \to 0.
  $$
 The second map is induced by restriction to $p$. The first map is
 multiplication by a holomorphic function vanishing once at $p$.
\end{lem}

\begin{proof}
 We have a defining exact sequence for the bundle of parabolic
 homomorphisms from $E$ to $F$ given by
  $$
   0\to\PH(E,F) \to \Hom (E,F) \to
   \bigoplus_{x\in D} \frac{\Hom(E_x,F_x)}{\PH(E_x,F_x)} \to 0.
  $$
 Now we tensor with the skyscraper sheaf $\CC(p)$, to get
  $$
   0\to \Tor \left(\frac{\Hom(E_p,F_p)}{\PH(E_p,F_p)},\CC(p)\right)
   \to \PH(E,F)_p \to \Hom (E,F)_p \to
   \frac{\Hom(E_p,F_p)}{\PH(E_p,F_p)} \to 0.
  $$
This is because $\Tor \left(\frac{\Hom(E_x,F_x)}{\PH(E_x,F_x)},\CC(p)\right)=0$
for $p\neq x$, and  the fact that if
$\Theta$ is a torsion sheaf supported scheme-theoretically at $p$ (i.e.,
 supported at $p$ and with no infinitesimal information), we have that
 $\Tor (\Theta,\CC(p))\cong \Theta \otimes \cO(-p)$ naturally
(to see this, tensor
 the exact sequence $\cO(-p) \to \cO\to \CC(p)$ with $\Theta$). Hence
  $$
   0\to \frac{\Hom(E_p,F_p)}{\PH(E_p,F_p)} \otimes \cO(-p) \to
    \PH(E,F)_p \to \Hom(E_p,F_p) \to
   \frac{\Hom(E_p,F_p)}{\PH(E_p,F_p)} \to 0,
  $$
which yields
  $$
   0 \to \frac{\Hom(E_p,F_p)}{\PH(E_p,F_p)} \otimes \cO(-p)
   \to \PH(E,F)_p \to \PH(E_p,F_p) \to 0.
  $$
Locally, with a local coordinate $z$ vanishing at $p$, the second
map is given by $(f_0+f_1 z +\cdots)_p \mapsto f_0$. The first map
is $f_1 \mapsto (f_1 z)_p$.
\end{proof}

To clarify the Lemma, let us see an example, where $E$ has rank
$3$ and weights $\beta_i$, $F$ has rank $4$ and weights
$\alpha_j$, and
$\beta_1<\alpha_1<\alpha_2<\alpha_3<\beta_2<\beta_3<\alpha_4$.
Then a typical parabolic homomorphism from $E$ to $F$ has matrix
of the form:
 $$
 \phi(z)= \left(\begin{array}{>{\columncolor[rgb]{0.81,0.81,0.81}}ccc}\cline{1-1}
 \multicolumn{1}{>{\columncolor[rgb]{0.81,0.81,0.81}}c|}{\phi_{11}(z)} &  \phi_{12}(z) & \phi_{13}(z) \\
 \multicolumn{1}{>{\columncolor[rgb]{0.81,0.81,0.81}}c|}{\phi_{21}(z)} &  \phi_{22}(z) & \phi_{23}(z) \\
 \multicolumn{1}{>{\columncolor[rgb]{0.81,0.81,0.81}}c|}{\phi_{31}(z)} &
 \multicolumn{1}{c}{\phi_{32}(z)} &  \multicolumn{1}{c}{\phi_{33}(z)} \\
 \cline{2-2} \cline{3-3}
 \rowcolor[rgb]{0.81,0.81,0.81}  \phi_{41}(z) &  \phi_{42}(z) &
 \multicolumn{1}{>{\columncolor[rgb]{0.81,0.81,0.81}}c|}{\phi_{43}(z)} \\
 \end{array}\right)
 $$
around $p$. The parabolicity of $\phi$ means that for $z=0$, the
only non-zero entries are those below the broken line. The line in
the matrix is easy to construct: starting by the upper-left
corner, draw a horizontal line for each $\beta_j$, and a vertical
line for each $\alpha_i$, considering the $\alpha$'s and $\beta$'s
in increasing order. The sheaf $\PH(E,F)$ is actually a bundle
(since it is torsion-free) of rank $\rk(E)\rk(F)$. Its stalk at
$p$, $\PH(E,F)_p$, is formed by the matrices with entries which
are complex numbers below the broken line, and which are complex
numbers times $z$ above the line.

\begin{prop} \label{prop:conditions-chi=0}
 Assume that $g>0$, $\sigma_c\geq 2g-2$ and $r_{1}> r_{2}$. Let $T'$,
 $T''$ be $\sigma_c^\pm$-stable triples
 with $\mu_{\sigma_c}(T')=\mu_{\sigma_c}(T'')$. Then
 $\chi(C^\bullet(T'',T'))=0$ if and only if the following conditions
 hold:
  \begin{enumerate}
   \item[(1)] $E_2'=0$. 
   \item[(2)] $\phi'':E_2'' \to E_1''(D)$ is a fibre bundle isomorphism at
   $X\setminus D$. In particular, $r_2''=r_1''$. 
   \item[(3)] At any point $p\in D$, write $\phi''=z^{-1}
   (\phi_0 +\phi_1 z + \phi_2 z^2 + \cdots)$,
   where $z$ is a local holomorphic coordinate around $p$ in $X$. Then
   $\PH(E''_{1,p},E'_{1,p}) \to \PH(E_{2,p}'',E'_{1,p})$,
   $f\mapsto -f\circ \phi_0$, is surjective. 
   \item[(4)] At any $p\in D$, consider the induced homomorphism
   $\phi_1: \ker \phi_0 \to \coker \phi_0$. Then
   $\PH(\coker\phi_0,E'_{1,p}) \to
   \Hom(\ker \phi_0,E'_{1,p})$,
   $f\mapsto -f\circ \phi_1$, is surjective.
  \end{enumerate}
\end{prop}

\begin{proof}
By Theorem \ref{thm:codim}, $\chi(C^\bullet(T'',T'))=0$ if and
only if $a_1$ is an isomorphism. By Lemma  \ref{lem:geniso}, if
$a_1$ is generically an isomorphism and $r_{1}> r_{2}$ then
$E'_2=0$. This proves (1). Also $\phi'':E_2''\to E_1''(D)$ is
generically an isomorphism. Moreover the two bundles involved in
the complex $C^\bullet (T'',T')$ must be of the same rank and
of the same degree. The complex $C^\bullet(T'',T')$ reduces to
 $$
 \PH(E''_{1},E'_{1})
 \stackrel{a_1}{\longrightarrow}
  \PH(E_{2}'',E'_{1}(D)),
 $$
where $a_1 ( f ) = -f \circ \phi''$ is an isomorphism of bundles.
Restricting $a_1$ to the open subset $U=X\setminus D$, we have
that $\Hom(E''_{1},E'_{1})|_U \to \Hom(E_{2}'',E'_{1}(D))|_U$ is
an isomorphism. Hence $E_2''|_U\to E_1''(D)|_U$ is an isomorphism
of bundles, and  (2) follows.

Now let $p\in D$, take a neighbourhood $U$ of $p$, and a coordinate
$z$ vanishing at $p$. Hence we may write $\phi''=\phi_0 z^{-1} +
\phi_1 + \phi_2 z + \cdots$, where $\phi_i\in \Hom (E_{2,p}'',
E''_{1,p})$ and $\phi_0 \in \PH(E_{2,p}'',E''_{1,p})$, on $U$. We
want to characterize when
 $$
 \PH(E''_{1},E'_{1})_p \to
 \PH(E_{2},E'_{1}(D))_p=\PH(E_2'',E_1'(p))_p
 $$
is an isomorphism of vector spaces. It is enough to analyze when
this map is surjective. Using Lemma \ref{lem:extra}, we have a
commutative diagram whose rows are short exact sequences:
  $$
  \begin{array}{ccccc}
   \displaystyle\frac{\Hom(E''_{1,p},E'_{1,p})}{\PH(E''_{1,p},E'_{1,p})} \otimes \cO(-p)
   &\stackrel{\cdot z}{\longrightarrow} & \PH(E_1'',E_1')_p &\to& \PH(E''_{1,p},E'_{1,p}) \\
    \downarrow b_0 &&\downarrow  b_1 &&\downarrow b_2 \\
    \displaystyle \frac{\Hom(E_{2,p}'',E'_{1,p})}{\PH(E_{2,p}'',E'_{1,p})}
   &\stackrel{\cdot z}{\longrightarrow} & \PH(E_2'',E_1'(p))_p &\to& \PH(E_{2,p}'',E'_{1,p})\otimes \cO(p).
 \end{array}
 $$
The middle vertical arrow is induced by $f\mapsto -f\circ \phi''$.
Thus the right vertical arrow is induced by $f_0 \mapsto
-(f_0\circ \phi_0 ) z^{-1}$. The left vertical arrow is thus given
by $f_1 \mapsto -(f_1\circ \phi_0) z^{-1}$.

We want to characterize the cases where the middle vertical arrow
is surjective. Using the long exact sequence produced by the snake
lemma, we see that $b_1$ being surjective is equivalent to $b_2$
being surjective and the connecting homomorphism $\ker b_2 \to
\coker b_0$ also being surjective. The condition that $b_2$ is
surjective is exactly (3).

For the remaining condition, we need to spell out the connecting
homomorphism. Take $f_0 \in\PH(E''_{1,p},E'_{1,p})$ lying in
 $$
 \ker b_2= \PH(E''_{1,p}/\phi_0(E_{2,p}'') ,E'_{1,p}).
 $$
Lift $f_0$ to a local section of $\PH(E_1'',E_1')$ on $U$, e.g.
taking $f(z)\equiv f_0$. Compose with $\phi''$ to get
$-(f\circ\phi_0 +f\circ \phi_1 z +\cdots) z^{-1}$. Recalling that
$f\circ \phi_0=0$, the leading term is
 $$
 -f_0\circ\phi_1 \in \coker b_0
 =\frac{\Hom(E_{2,p}'',E'_{1,p})}{\PH(E_{2,p}'',E'_{1,p}) +
 b_0(\Hom(E''_{1,p},E'_{1,p}))}.
 $$
Assuming that (3) holds already, we have that
$\PH(E_{2,p}'',E'_{1,p}) \subset
b_0(\PH(E_{1,p}'',E'_{1,p}))\subset
b_0(\Hom(E_{1,p}'',E'_{1,p}))$, since the maps $b_0$ and $b_2$ are
both composition with $\phi_0$. Hence the image of $f_0$ under the
connecting homomorphism is 
 $$
 -f_0\circ\phi_1 \in \coker b_0 =
 \frac{\Hom(E_{2,p}'',E'_{1,p})}{b_0(\Hom(E''_{1,p},E'_{1,p}))}= \Hom(\ker
 \phi_0,E'_{1,p})\, .
 $$
Therefore the surjectivity of the connecting homomorphism is
equivalent to (4).
\end{proof}


\begin{lem} \label{lem:(4)}
 Condition {\rm (4)} of Proposition \ref{prop:conditions-chi=0}
 holds  if and only if all the weights of
 $E'_{1,p}$ are bigger than those of $\coker \phi_0$, and
 $\phi_1:\ker \phi_0 \to \coker\phi_0$ is an isomorphism.
\end{lem}

\begin{proof}
The condition (4) says that
 $$
 \PH \left(\frac{E_{1,p}''}{\phi_0(E_{2,p}'')}, E'_{1,p} \right) \to
  \Hom(\ker \phi_0, E'_{1,p}), \quad f\mapsto -f\circ \phi_1,
 $$
is surjective. Since $E_{1,p}''/\phi_0(E_{2,p}'')$ and
$\ker\phi_0$ are vector spaces of the same dimension, this is
equivalent to the following two conditions:
 \begin{itemize}
  \item $\phi_1: E_{2,p}''\to E_{1,p}'' \subset E_{1,p}$ satisfies
  that $\phi_1: \ker \phi_0 \to \coker\phi_0$ is an isomorphism.
  \item $\PH(E_{1,p}''/\phi_0(E_{2,p}''), E'_{1,p})=
  \Hom( E_{1,p}''/\phi_0(E_{2,p}''), E'_{1,p})$.
  Hence all the weights of $E''_{1,p}/\phi_0(E_{2,p}'')$ are smaller than
  those of $E_{1,p}'$. 
 \end{itemize}
\end{proof}

Let $\sigma_c\in (\sigma_m,\sigma_M)$ be a critical value with
$\sigma_c\geq 2g-2$. We aim to characterize when
$\cN_{\sigma_c^-}^s$ and $\cN_{\sigma_c^+}^s$ are birational by
using Proposition \ref{prop:codim}. Let us deal with either of
$\cS_{\sigma_c^\pm}$. Suppose that $T'$ and $T''$ are
$\sigma_c$-semistable, $\sigma_c^\pm$-stable triples with
$\mu_{\sigma_c}(T')=\mu_{\sigma_c}(T'')$. We consider extensions
 \begin{equation}\label{eqn:harto4}
 0 \to T''\to T\to T'\to 0
 \end{equation}
(note that we have changed the role of $T'$ and $T''$ in the
computation of the codimension of the flip loci in Section
\ref{sec:chi}, so that now $T''$ is the subtriple), where
$\mu_{\sigma_c^\pm}(T'')<\mu_{\sigma_c^\pm}(T)$, by Proposition
\ref{prop:triplesext}. The first conclusion to infer from
Proposition \ref{prop:conditions-chi=0} is that, if
$\chi(C^\bullet(T'',T'))=0$ then $r_2'=0$ and $r_2''=r_1''$. So
$\mu_{\sigma_c^+}(T'')>\mu_{\sigma_c^+}(T)$. Therefore
$\cS_{\sigma_c^+}$ cannot be of zero codimension. So our study is
limited to $\cS_{\sigma_c^-}$: the only situation we may encounter
when $\chi(C^\bullet(T'',T'))=0$ is that $\cN_{\sigma_c^-}^s$ has
more irreducible components than $\cN_{\sigma_c^+}^s$.

\medskip

To analyze when $\chi(C^\bullet(T'',T'))=0$ we have to check when
conditions (3) and (4) of Proposition \ref{prop:conditions-chi=0}
are satisfied. Let $p\in D$ be a parabolic point. We need to
understand the {\em parabolic vector spaces} $E_{2,p}$ and
$E_{1,p}$. These have parabolic weights of multiplicity one and
all weights are different, by Assumption \ref{assumption2}. We
shall keep the following notation for the rest of the section:
$\alpha_i$ denote the weights of $E_{1,p}$ and $\beta_j$ denote
the weights of $E_{2,p}$ (we drop $p$ from the notation in the
weights when this causes no confusion).

Since $T$ is a triple which is an extension (\ref{eqn:harto4})
with $r_2'=0$ and $r_2''=r_1''$, then $\phi:E_2\to E_1(D)$ comes
from a map $\phi'':E_2\to E_1''(D)$ as follows
 $$
 \begin{array}{ccccc}
  \,\,\,\,\, E_2 &=& \,\, E_2 & \to & 0\\
   \phi''\downarrow &&\phi\downarrow &&\downarrow \\
  E_1''(D) &\to & E_1(D) & \to & E_1'(D) .
 \end{array}
 $$
Take a neighbourhood $U$ of $p$ where $E_1|_U=E_1'|_U\oplus
E_1''|_U$. Then $\phi=(\phi_0+\phi_1 z+\cdots) z^{-1}$ and
$\phi_0:E_{2,p} \to E_{1,p}''$ is a parabolic map. This gives
decompositions of the parabolic vector spaces
  \begin{equation}\label{eqn:decomp}
  \begin{array}{l}
  E_{1,p}=E'_{1,p}\oplus E''_{1,p}, \\
  E''_{1,p}= \im \phi_0 \oplus \coker \phi_0,
  \end{array}
  \end{equation}
as direct sums of parabolic vector subspaces (the splitting is
non-canonical, but the weights of the different subspaces are
well-determined).

Let us see that there is a ``canonical'' distribution of weights
in (\ref{eqn:decomp}) such that conditions (3) and (4) hold. Note
that $\PH(E_{2,p},E_{1,p})$ is a vector space, in particular an
irreducible affine variety. We may consider the action of
$\ParAut(E_{2,p}) \times \ParAut(E_{1,p})$ on this space (this
corresponds to lower triangular changes of bases). Then there is a
unique open dense orbit, which is the only orbit of maximal
dimension. We shall call an element of such orbit a \emph{generic}
parabolic homomorphism of $E_{2,p}$ to $E_{1,p}$.  For instance,
if $E_{2,p}$ is $7$-dimensional with weights $\beta_j$ and
$E_{1,p}$ is $9$-dimensional with weights $\alpha_i$, and
 $$
 \alpha_1 <\beta_1 <\beta_2<\alpha_2 <\beta_3
 <\beta_4 <\alpha_3 <\alpha_4 <\alpha_5<\alpha_6 <\beta_5 <\alpha_7
 <\alpha_8 <\beta_6 <\beta_7 <\alpha_9,
 $$
then the generic elements are the orbit of the element
 \begin{equation}\label{eqn:example}
 \left( \begin{array}{ccccccc}
 \multicolumn{1}{|c}{0} &0 &0&0&0&0&0 \\
 \cline{1-2} 0&\multicolumn{1}{c|}{1}& 0&0&0&0&0 \\
 \cline{3-4} 0&0&0&\multicolumn{1}{c|}{1}&0&0&0 \\
  0&0&1&\multicolumn{1}{c|}{0}&0&0 &0\\
 1&0&0&\multicolumn{1}{c|}{0}&0&0 &0\\
 0&0&0&\multicolumn{1}{c|}{0}&0&0&0 \\
  \cline{5-5} 0&0&0&0&\multicolumn{1}{c|}{1}&0&0\\
 0&0&0&0&\multicolumn{1}{c|}{0}&0 &0 \\
  \cline{6-7} 0&0&0&0&0 &0&\multicolumn{1}{c|}{1}\end{array} \right).
 \end{equation}

\begin{lem} \label{lem:phi0generic}
  Suppose that $\phi_0:E_{2,p}\to E_{1,p}$ is a generic parabolic
  homomorphism, and let $E_{1,p}=E_{1,p}'\oplus
  E_{1,p}''$ be any parabolic splitting with
  $\im\phi_0 \subset E_{1,p}''$. Then condition {\rm (3)} in
  Proposition \ref{prop:conditions-chi=0} is  satisfied.
\end{lem}

\begin{proof}
Suppose that $\phi_0$ is a generic element in
$\PH(E_{2,p},E_{1,p})$, and let us see that  the map
$\PH(E''_{1,p},E'_{1,p}) \to \PH(E_{2,p},E'_{1,p})$,  $f\mapsto
-f\circ \phi_0$, is surjective. Take $g\in \PH(E_{2,p},E'_{1,p})$.
Consider the map $\phi_\epsilon=\phi_0 \oplus \epsilon g:
E_{2,p}\to E''_{1,p}\oplus E'_{1,p}$. For $\epsilon$ small we have
that $\phi_\epsilon$ also lives in the generic open set, so it is
equivalent to $\phi_0$ by the action of $\ParAut(E_{2,p}) \times
\ParAut(E_{1,p})$. This means that
  $$
  \left(\begin{array}{cc} a_\epsilon &b_\epsilon \\
  c_\epsilon & d_\epsilon\end{array}\right)
  \left(\begin{array}{c} \phi_0 \\
  0\end{array}\right) M_\epsilon =
  \left(\begin{array}{c} \phi_0 \\ \epsilon g
  \end{array}\right).
  $$
Both matrices, $\left(\begin{array}{cc} a_\epsilon &b_\epsilon \\
c_\epsilon & d_\epsilon\end{array}\right)$ and $M_\epsilon$, are
the identity for $\epsilon =0$, so $a_\epsilon$ is invertible for
small $\epsilon$. Therefore $\phi_0 M_\epsilon=a_\epsilon^{-1}
\phi_0$ and $c_\epsilon\phi_0 M_\epsilon= \epsilon g$. This yields
 $$
 g= \epsilon^{-1} c_\epsilon a_\epsilon^{-1} \phi_0,
 $$
as required.
\end{proof}

Recall that we have fixed topological data (fixed ranks, degrees
and parabolic weights) for the triples $T$ we are studying. When
we write such a triple $T$ as an extension $T''\to T \to T'$,
there are different possible topological types for $T'$ and $T''$.
By the above discussion, our best chance to obtain
$\chi(C^\bullet(T',T''))=0$ is to arrange the topological types as
follows:
 \begin{itemize}
 \item Fix the ranks $r_2'=0$, $r_2''=r_2$, $r_1''=r_2$,
 $r_1'=r_1-r_2$. This is necessary for conditions (1) and (2) to
 hold. So $\phi:E_2\to E_1(D)$ should be induced by
 $\phi'':E_2''\to E_1''(D)$ by means of the inclusion $E_1''(D)\to
 E_1(D)$.
 \item At each $p\in D$, consider a generic element $\phi_p\in
 \PH(E_{2,p},E_{1,p})$. This determines the weights of $\im
 \phi_p \subset E_{1,p}''$. By Lemma \ref{lem:phi0generic}
 condition (3) is satisfied.
 \item Choose the weights of $\coker \phi''_p$ in the unique way
 such that Lemma \ref{lem:(4)} is satisfied. This gives the
 weights of $E''_{1,p}=\im \phi_p \oplus \coker \phi''_p$ at each
 $p\in D$, and hence the weights of $E_{1,p}'$.
 \item $d_2''=d_2$. Now condition (2) determines the degree of $E_1''$,
 since the map $\phi'':E_2\to E_1''(D)$ is an isomorphism on $X\setminus D$ and
 it is of a specified form at each $p\in D$. Namely, introduce the
 number
  \begin{equation}\label{eqn:rp}
  r_p = \min \{ \dim \coker \psi_0 \ | \ \psi_0\in
  \PH(E_{2,p},E_{1,p})\} - (r_1-r_2).
  \end{equation}
 Obviously this minimum is obtained for a generic parabolic morphism.
 Moreover $r_p=\dim \coker \phi_0$, where $\phi_p:E_{2,p}\to
 E_{1,p}$ is generic, and $\phi_0=\phi_p:E_{2,p}\to E_{1,p}''$,
 using that $E_{1,p}''\subset E_{1,p}$.
 With this notation, $E_2\to E_1''(D) \to \oplus_{p\in D} \CC(p)^{r_p}$ is
 an exact sequence of sheaves, so $d_1''=d_2-r_2s+\sum_{p\in D} r_p$.
 \end{itemize}

This does not guarantee the existence or uniqueness of the
topological types of $T'$ and $T''$ to have
$\chi(C^\bullet(T',T''))=0$, but helps us in which direction to
look for such distributions of topological types.

Let us see this discussion in the particular example
(\ref{eqn:example}). For a generic $\phi_p:E_{2,p}\to E_{1,p}$,
the weights of $\im \phi_0$ are
$\alpha_2,\alpha_3,\alpha_4,\alpha_5,\alpha_7,\alpha_9$, and the
weight of $\coker\phi_0$ is $\alpha_1$. Thus the weights of
$E'_{1,p}$ are $\alpha_6,\alpha_8$. The map $\phi$ takes the form:
 $$
 \phi= \Big( \, \left( \begin{array}{ccccccc}
 \multicolumn{1}{|c}{0} &0 &0&0&0&z&0 \\
 \cline{1-2} 0&\multicolumn{1}{c|}{1}& 0&0&0&0&0 \\
 \cline{3-4} 0&0&0&\multicolumn{1}{c|}{1}&0&0&0 \\
  0&0&1&\multicolumn{1}{c|}{0}&0&0 &0\\
 1&0&0&\multicolumn{1}{c|}{0}&0&0 &0\\
 0&0&0&\multicolumn{1}{c|}{0}&0&0&0 \\
  \cline{5-5} 0&0&0&0&\multicolumn{1}{c|}{1}&0&0\\
 0&0&0&0&\multicolumn{1}{c|}{0}&0 &0 \\
  \cline{6-7} 0&0&0&0&0 &0&\multicolumn{1}{c|}{1}\end{array} \right) + O(z^2) \,
  \Big) \cdot z^{-1},
 $$
around $p\in D$. Note that such $\phi:E_2\to E_1(D)$ is injective
for $z\neq 0$, as required by condition (2).

\begin{rmk}\label{rem:rp}
 The definitions of generic parabolic map and of $r_p$ given in
 (\ref{eqn:rp}) are also valid in the case $r_1=r_2$.
\end{rmk}

\begin{prop} \label{prop:atmost}
 Assume $g>0$, $r_1>r_2$ and $\sigma_c\geq 2g-2$. Let $T'$, $T''$ be
 $\sigma_c^-$-stable triples
 with $\mu_{\sigma_c}(T')=\mu_{\sigma_c}(T'')$.
 If $\chi(C^\bullet(T'',T'))=0$ then the following holds:
 \begin{enumerate}
   \item[(i)] $r_2'=0$, $r_2''=r_2=r_1''$, $d_2''=d_2$.
   \item[(ii)] For each $p\in D$, the parabolic map $\phi''_p:E_{2,p}\to
   E''_{1,p}$ has rank $r_2-r_p$, with $r_p$ defined in
   (\ref{eqn:rp}).
   \item[(iii)] $d_1''=d_2-r_2s+\sum_{p\in D} r_p$.
   \item[(iv)] For each $1\leq k\leq r_2-r_p$, define
     \begin{equation}\label{eqn:ik}
     i_k=\min\{ j  \ |\ 1\leq j\leq r_1, \beta_k <\alpha_j, j>i_{k-1}\}.
     \end{equation}
   and let $I=\{i_1,\ldots, i_{r_2-r_p}\}$. Let $J\subset
   \{1,\ldots, r_1\} - I$ be the set of the lowest $r_p$ elements
   of $\{1,\ldots, r_1\} - I$. Then the weights of $E''_{1,p}$ are
   exactly $\{\alpha_i \ | \ i\in I\cup J\}$.
 \end{enumerate}
  In particular, the ranks, degrees and weights of $T'$ and $T''$
  are univocally determined. Thus there is at most one possible
  value of $\sigma_c$ for which $\chi(C^\bullet(T'',T'))=0$.
\end{prop}

\begin{proof}
Item (i) follows from Proposition \ref{prop:conditions-chi=0} (1).

Item (iii) follows once we know item (ii) and using Proposition
\ref{prop:conditions-chi=0} (2), since in this case we have an
exact sequence of sheaves
  $$
   E_2 =E_2''\stackrel{\phi''}{\longrightarrow} E_1''(D) \to \bigoplus_{p\in
   D} \CC(p)^{r_p}.
  $$

Next, note that the increasing sequence of numbers $i_1,i_2,\ldots
\in \{1,\ldots, r_1\}$ is well-defined for $1\leq k\leq r_2-r_p$.
Actually, looking at a generic parabolic map $\psi_0: E_{2,p}\to
E_{1,p}$, the weights of $\im \psi_0$ are
$\alpha_{i_1},\ldots,\alpha_{i_{r_2-r_p}}$, with $r_2-r_p=\dim \im
\psi_0$ (see (\ref{eqn:example}) for a specific example).

Now we shall prove (ii) and (iv) using Proposition
\ref{prop:conditions-chi=0} (3), i.e., that
  \begin{eqnarray}\label{eqn:surj}
   \PH(E''_{1,p},E'_{1,p}) &\to &\PH(E_{2,p},E'_{1,p}) \notag\\
    f &\mapsto & f\circ \phi_0
  \end{eqnarray}
is surjective, denoting as before, $\phi_0=\phi''_p$. Let
$\{e_1,\ldots, e_{r_1}\}$ be a basis for $E_{1,p}$ adapted to its
parabolic structure (and adapted to the splitting $E'_{1,p}\oplus
E''_{1,p}$, i.e. each $e_i$ belongs either to $E'_{1,p}$ or
$E''_{1,p}$), and let $\{v_1,\ldots, v_{r_2}\}$ be a basis for
$E_{2,p}$ adapted to its parabolic structure.

Now let $t_0\in \{1,\ldots, r_1\}$ such that $\alpha_{t_0}$ is the
lowest weight of $E'_{1,p}$. Let $0\leq a\leq r_2-r_p$ such that
$i_a<t_0\leq i_{a+1}$ (introducing the notation $i_0=0,
i_{r_2-r_p+1}=r_1+1$). Let us see that $\alpha_{i_{a+1}}, \ldots,
\alpha_{r_2-r_p}$ are weights of $\im\phi_0$ (if $a=r_2-r_p$ then
there is nothing to prove). Actually, they cannot be weights of
$\coker\phi_0$, since by Lemma \ref{lem:(4)} all the weights of
$\coker \phi_0$ are smaller than $\alpha_{t_0}$. So they are
weights of $\im \phi_0$ or of $E_{1,p}'$ by (\ref{eqn:decomp}).
Suppose that $\alpha_{i_{a+1}}, \ldots, \alpha_{i_{b-1}}$ are
weights of $\im\phi_0$ but $\alpha_{i_b}$ is the first weight of
$E'_{1,p}$ in the list. Then take $V=\langle v_1,\ldots, v_{b}
\rangle \subset E_{2,p}$. The surjectivity of (\ref{eqn:surj})
gives that
  $$
  \PH(E''_{1,p},\langle e_{i_b}\rangle) \twoheadrightarrow
  \PH(E_{2,p},\langle e_{i_b}\rangle) \twoheadrightarrow
  \PH(V,\langle e_{i_b}\rangle) =\Hom(V,\langle e_{i_b}\rangle)
  $$
is surjective (the last equality follows from
$\alpha_{i_b}>\beta_b$). Therefore $\phi_0|_V:V\to E''_{1,p}$ must
be injective, and all the weights or $\phi_0(V) \subset E''_{1,p}$
should be smaller than $\alpha_{i_b}$. So there are weights
$\alpha_{x_1}<\ldots <\alpha_{x_{b}}<\alpha_{i_b}$ with
$\beta_j<\alpha_{x_j}$. This implies that $i_j \leq x_j$,
$j=1,\ldots, b$, which contradicts that $x_b<i_b$.

Next step is to see that there are $y_1<\cdots <y_a<t_0$ such that
$i_j\leq y_j$, $j=1,\ldots, a$ and $\alpha_{y_j}$ are weights of
$\im\phi_0$. As before, take $V=\langle v_1,\ldots, v_{a} \rangle
\subset E_{2,p}$. The surjectivity of (\ref{eqn:surj}) gives that
  $$
  \PH(E''_{1,p},\langle e_{t_0}\rangle) \twoheadrightarrow
  \PH(E_{2,p},\langle e_{t_0}\rangle) \twoheadrightarrow
  \PH(V,\langle e_{t_0}\rangle) =\Hom(V,\langle e_{t_0}\rangle)
  $$
is surjective. So $\phi_0|_V:V\to E''_{1,p}$ must be injective,
and all the weights or $\phi_0(V) \subset E''_{1,p}$ should be
smaller than $\alpha_{t_0}$. So there are weights
$\alpha_{y_1}<\ldots <\alpha_{y_{a}}<\alpha_{t_0}$ with
$\beta_j<\alpha_{y_j}$. This implies that $i_j \leq y_j$,
$j=1,\ldots, a$.

The elements
  \begin{equation}\label{eqn:weight-im}
  \{ y_1,\ldots, y_a, i_{a+1}, \ldots, i_{r_2-r_p}\}
  \end{equation}
are weights of $\im \phi_0$. So $\dim\im \phi_0\geq r_2-r_p$. As
obviously  $\dim\im \phi_0\leq r_2-r_p$, it must be $\dim\im
\phi_0 =r_2-r_p$, implying item (ii). Thus the weights of
$\im\phi_0$ are exactly those in (\ref{eqn:weight-im}). The
elements
  \begin{equation} \label{eqn:weight-coker}
  \{1,\ldots, t_0-1\} - \{y_1,\ldots, y_a\}
  \end{equation}
are the sub-indices of the weights of $\coker\phi_0$, by Lemma
\ref{lem:(4)}. So $t_0-1-a=r_p$, i.e. $t_0=r_p+a+1$. Finally (the
sub-indices of) the weights of $E''_{1,p}$ are
  \begin{eqnarray*}
   &&(\{1,\ldots, t_0-1\} - \{y_1,\ldots, y_a\}) \cup \{ y_1,\ldots, y_a, i_{a+1}, \ldots,
  i_{r_2-r_p}\}= \\
  &&=\{1,\ldots, t_0-1\} \cup \{i_{a+1}, \ldots,
  i_{r_2-r_p}\} = I \cup J,
  \end{eqnarray*}
as required.
\end{proof}

Our final result in this section completes the proof of Theorem
\ref{thm:connected-minima}. We have to use Theorem
\ref{thm:r1=r2}, which will be proven in the next section. First,
consider the distribution of weights and degrees given by
Proposition \ref{prop:atmost}, and consider the critical value
associated to it, which is
 \begin{equation}\label{eqn:sL}
 \begin{aligned}
  \sigma_L &=
 \left(1+\frac{r_{1}+r_{2}}{r_{1}-r_{2}}\right)
 (\pmu(E_{1})-\pmu(E_{2}))+s
 \frac{r_{1}+r_{2}}{r_{1}-r_{2}} -\frac{\pdeg (E_1''(D))-\pdeg
 (E_2)}{r_2}  \\ &=\sigma_M- \frac{1}{r_2}\, \epsilon\, ,
 \end{aligned}
 \end{equation}
where
 \begin{equation}\label{eqn:epsilon}
  \epsilon= \pdeg (E_1''(D))-\pdeg (E_2) >0,
 \end{equation}
and the weights and degree of $E_1''$ are given by Proposition
\ref{prop:atmost}. For instance, in the example worked out in
(\ref{eqn:example}), $\epsilon= \sum_{i\neq 6,8} \alpha_i -\sum
\beta_j +1$.

The value of $\sigma_L$ is very close to $\sigma_M$ but strictly
smaller, as expected.

\begin{thm}\label{thm:mainharto}
  Assume $r_1>r_2$ and $g>0$. If $\sigma_L>2g-2$ then
  $\cN_\sigma^s$ is irreducible and non-empty for all
  $2g-2\leq\sigma <  \sigma_L$. If $\sigma_L<2g-2$ then
  $\cN_\sigma^s$ is empty for all $\sigma\ge 2g-2$.
\end{thm}

\begin{proof}
First, note that for $\sigma>\sigma_M$, $\cN_{\sigma}$ is empty by
Proposition \ref{prop:cotasigma}. Assume for a while that
$\cN_{\sigma}^s$ is non-empty for some value of $\sigma\geq 2g-2$,
then there must exist the minimum value $\tilde\sigma_L\in
(2g-2,\sigma_M)$ of $\sigma$ such that
$\cN^s_{\tilde\sigma_L^+}=\emptyset$ and
$\cN^s_{\tilde\sigma_L^-}\neq \emptyset$. Clearly this
$\tilde\sigma_L$ is a critical value and it must correspond to a
set of extensions $T''\to T\to T'$ with
$\chi(C^\bullet(T'',T'))=0$.

By Proposition \ref{prop:atmost} there is at most one
(topological) possibility for $T'$ and $T''$ to have
$\chi(C^\bullet(T'',T'))=0$. This implies that
$\tilde\sigma_L=\sigma_L$. For any other critical value
$\sigma_c$, the moduli spaces $\cN_{\sigma_c^+}^s$ and
$\cN_{\sigma_c^-}^s$ are birational, by Proposition
\ref{prop:codim}. So all moduli spaces $\cN_\sigma^s$ are
birational for $2g-2\leq \sigma<\sigma_L$.

Moreover there may be different distributions of weights, ranks
and degrees giving rise to the critical value $\sigma_L$, but only
the one given by Proposition \ref{prop:atmost} gives critical
subsets $\cS_{\sigma_L^-}$ of codimension zero. So the number of
irreducible components is given by the number of irreducible
components of a subset of the space of extensions $T''\to T\to T'$
with the distribution of weights, ranks and degrees given by
Proposition \ref{prop:atmost}. Let us see that this space of
extensions is non-empty and irreducible: the triples $T'$ have
$r_2'=0$, $r_1'=r_1-r_2$, so they are parametrized by a moduli
space of parabolic bundles $E_1'$, which is non-empty, irreducible
and of the expected dimension by \cite{by2}. The triples $T''$
have $r_1''=r_2''=r_2$, and $d_1''+r_1''s-d_2''-\sum r_p=0$, so
they are parametrized by a moduli space of $\sigma_L^-$-stable
triples which is non-empty, irreducible and of the expected
dimension by Theorem \ref{thm:r1=r2}. Now the dimension of the
projective fibres of the space of extensions $T''\to T\to T'$ is
 $$
 -\chi(C^\bullet(T',T''))-1 \geq 0,
 $$
since $\chi(C^\bullet(T',T''))<0$, by Theorem \ref{thm:codim}.
Therefore there is a non-empty space of extensions. Moreover,  a
generic triple $T'$ is $\sigma_L$-stable. In that case, any
non-trivial extension $T''\to T\to T'$ is $\sigma_L^-$-stable (see
Proposition \ref{prop:triplesext}). So the space
$\cS_{\sigma_L^-}$ is non-empty, and irreducible.

Finally, if $\sigma_L>2g-2$, the argument above proves that
$\cN_{\sigma_L^-}^s$ is non-empty, so there is some non-empty
$\cN_\sigma$ with $\sigma>2g-2$ and the statement of the theorem
follows. Conversely, if some $\cN_\sigma$ with $\sigma>2g-2$ is
non-empty, then it must be $\sigma_L>2g-2$ completing the
argument.
\end{proof}

Now Proposition \ref{prop:harto11} transfers the inequalities
$\sigma_m\leq 2g-2<\sigma_L$ into a Milnor--Wood type inequality
$0\leq |\tau|<\tau_L$, where
 \begin{equation}\label{eqn:tauL}
  \tau_L= \min\{p,q\}(2g-2+s) - \frac{|p-q|}{p+q} \epsilon\, ,
 \end{equation}
where $\epsilon$ is given in \eqref{eqn:epsilon}.

\begin{rmk}\label{rem:hartoxx}
One can spell out the process for computing $\epsilon$, by using
the procedure of Proposition \ref{prop:atmost} and the
identification of Proposition \ref{prop:iden}. Let $p=\rk(V)$,
$q=\rk(W)$, $\alpha$ the system of weights of $V$ and $\beta$ the
system of weights of $W$. Suppose that $q \leq p$ (the other case
is similar, interchanging the roles of $V$ and $W$). Define, at
each $x\in D$, $\alpha_{i+pl}(x)=\alpha_i(x)+l$, for any $l\geq
1$. Put $i_0=0$ and define, for $1\leq k\leq q$,
  $$
  i_k=\min \{ j \ | \ j>i_{k-1}, \alpha_j >\beta_k \}.
  $$
Then
  $$
  \epsilon= \sum_{x\in D}\sum_{k=1}^p  (
  \alpha_{i_k}(x)-\beta_k(x)) .
  $$
\end{rmk}

\section{The moduli space of triples for $r_1=r_2$ and large $\sigma$}
\label{sec:r_1=r_2}

In this section, we study the moduli space of triples with equal
ranks $r_1=r_2$. We prove that some of them are irreducible and
non-empty for $\sigma\geq 2g-2$. The results here are enough for
the proof of Theorem \ref{thm:mainharto} to work, but we also
analyze some other cases. It is likely that the result holds in
general.

\begin{prop}\label{prop:equalrs}
  Suppose that $r_1=r_2$ and $g>0$. Then all the moduli spaces $\cN_\sigma$, for
  $\sigma\geq 2g-2$ are birational to each other.
\end{prop}

\begin{proof}
 This is a consequence of Theorem \ref{thm:codim} and Proposition
 \ref{prop:codim}. For
 $\chi(C^\bullet(T',T''))$ to vanish, it must be $a_1$ an
 isomorphism. But this is impossible if $r_1=r_2$ by Lemma
 \ref{lem:geniso}.
\end{proof}

Now let us see that the moduli spaces $\cN_\sigma$ stabilizes for
$\sigma$ large.

\begin{prop} \label{prop:inj}
Suppose that $r_1=r_2$. Then there is a value $\sigma_1$ such that
any $\sigma$-stable parabolic triple $T=(E_1,E_2,\phi)$ with
$\sigma>\sigma_1$ satisfies that $\phi$ is injective. Hence
 \begin{equation}\label{eqn:harto12}
 0\to E_{2}\to E_{1}(D)\to S\to 0,
 \end{equation}
where $S$ is a torsion sheaf.
\end{prop}

\begin{proof}
Denote $N=\ker\phi$ and consider the parabolic subtriple
$(0,N,\phi)$. Suppose that $k=\rk(N)>0$. The $\sigma$-stability of
$T$ implies that
  $$
  \pdeg N + k\sigma < k
  \left( \frac{\pdeg(E_{1}\oplus E_{2})}{2r_1}
  +\frac12 \sigma \right).
  $$
Now consider the subtriple $(I, E_{2},\phi)$ where $I(D)$ is the
parabolic image sheaf of $\phi$, with rank $\rk(I)=r_1-k$. The
$\sigma$-stability of $T$ gives us
 $$
 \pdeg (I\oplus E_2) + r_1 \sigma < (2r_1- k)\left(
 \frac{\pdeg(E_{1}\oplus E_{2})}{2r_1}
 +\frac12 \sigma \right).
 $$
Adding up both equations, and noting that $\pdeg N +\pdeg
I(D)=\pdeg E_2$, we get
 $$
 2 \pdeg E_2 -(r_1-k)s +(r_1+k)\sigma < \pdeg(E_{1}\oplus E_{2}) +
 r_1\sigma,
 $$
which is rewritten as
 $$
 \sigma \leq \frac{\pdeg E_1-\pdeg E_2 +
 (r_1-k)s}{k}\ .
 $$
So for $\sigma_1=\pdeg E_1-\pdeg E_2 + (r_1-1)s$ the result
follows.
\end{proof}


\begin{lem} \label{lem:t'}
 Suppose that $r_1=r_2$ and $\sigma>\sigma_1$. Let $T$ be a
 $\sigma$-stable triple and $T'$ a subtriple of $T$ with $r_1'=r_2'$.
 Write $E_2\to E_1(D)\to S$,  $E_2'\to E_1'(D)\to S'$, $t=\length
 S$, $t'=\length S'$. Then
 \begin{eqnarray*}
  \pmu(E'_{1})&<&
  \pmu(E_{1})+\frac{1}{2}\left(\frac{t'}{r'_1}-\frac{t}{r_1}\right)+s ,\\
  \pmu(E'_{2})&<&
  \pmu(E_{2})-\frac{1}{2}\left(\frac{t'}{r'_1}-\frac{t}{r_1}\right)+s
  .
 \end{eqnarray*}
\end{lem}

\begin{proof}
{}From Proposition \ref{prop:inj}, as $\sigma>\sigma_{1}$, $\phi$
is an injective morphism. So $\phi'$ is injective for any
subtriple $T'$ of $T$. Hence for a subtriple $T'$ with
$r'_{1}=r'_{2}$ we have the following commutative diagram
 \begin{displaymath}
 \begin{CD}
  0 @>>> E'_{2}@>>> E'_{1}(D) @>>> S' @>>> 0\\
  @.   @VVV @VVV @VVV\\
  0@>>> E_{2}@>>> E_{1}(D) @>>> S @>>> 0,
\end{CD}
\end{displaymath}
where $S$ and $S'$ are torsion sheaves. Let  $t$ and $t'$ denote
the lengths of $S$ and $S'$ respectively, as in the statement. By
stability,
\begin{eqnarray*}
 0 & > & \pmu_{\sigma}(T')-\pmu_{\sigma}(T) \\
 &=& \frac12 \left(\pmu(E'_{1})+\pmu(E'_{2})-\pmu(E_{1})-\pmu(E_{2})\right) \\
 &=& \pmu(E'_{1})-\pmu(E_{1})-\frac{1}{2}(\pmu(E'_{1}) -
 \pmu(E'_{2})) + \frac12 (\pmu(E_{1}) - \pmu(E_{2})) \\
 &=& \pmu(E'_{2})-\pmu(E_{2})+\frac{1}{2}(\pmu(E'_{1}) -
 \pmu(E'_{2})) - \frac12 (\pmu(E_{1}) - \pmu(E_{2})).
\end{eqnarray*}
Now at each point $p\in D$, $|\sum \beta_j(p) -\sum \alpha_i(p)|
\leq r_1$, so $t-r_1 s \leq \pdeg E_1(D) -\pdeg E_2 \leq t+ r_1
s$, equivalently $t-2r_1 s \leq \pdeg E_1-\pdeg E_2 \leq t$ or
 $$
 \frac{t}{r_1} -2 s \leq \pmu( E_1)-\pmu( E_2) \leq \frac{t}{r_1}.
 $$
Analogously, for $T'$ we have
 $$
 \frac{t}{r_1'} -2 s \leq \pmu( E_1')-\pmu( E_2') \leq \frac{t}{r_1'}.
 $$
Substituting into the formulae above, we get the result in the
statement.
\end{proof}

\begin{prop} \label{prop:stabilize}
Suppose that $r_1=r_2$. Then there is a value $\sigma_2\geq
\sigma_1$ such that $\cN_{\sigma}^s=\cN_{\sigma'}^s$ for any
$\sigma,\sigma'\geq \sigma_2$, i.e. there are no critical values
above $\sigma_2$.
\end{prop}

\begin{proof}
Consider a $\sigma$-stable triple $T=(E_1,E_2,\phi)$ with
$\sigma>\sigma_{1}$. Suppose that $T$ is properly
$\sigma_c$-semistable for some $\sigma_c$, and let $T'\subset T$
be a $\sigma_c$-destabilizing subtriple. Clearly $r_2'\leq r_1'$,
since $\phi$ being injective implies that $\phi'$ is also
injective. On the other hand, if $r_1'=r_2'$ then  $T$ is
$\sigma$-semistable for generic values of $\sigma$ and could not
be $\sigma$-stable for some $\sigma$. Therefore $r_2'<r_1'$. In
the formula
 \begin{equation}\label{eqn:sigma}
 \sigma_{c}=2
 \pmu(E'_{1})\frac{r'_{1}}{r'_{1}-r'_{2}}+2\pmu(E'_{2})\frac{r'_{2}}{r'_{1}-r'_{2}}
 -(\pmu(E_{1})+\pmu(E_{2}))\frac{r'_{1}+r'_{2}}{r'_{1}-r'_{2}},
\end{equation}
we want to bound the values of $\pmu(E'_{1})$ and $\pmu(E'_{2})$
in order to get a bound for the critical value $\sigma_{c}$ which
is independent of $T$.

Apply Lemma \ref{lem:t'} to the subtriples
$(\phi'(E_{2}')(-D),E'_{2},\phi')$ and $(E_1',
(\phi')^{-1}(E_{1}'(D)),\phi')$, both of which satisfy the equal
rank condition. The first one has no torsion, the second has
torsion with $0\leq t'\leq t$. We get
 \begin{eqnarray*}
  \pmu(E'_{2})&<& \pmu(E_{2})+ \frac{t}{2r_1}+s, \\
 \pmu(E'_{1})&<& \pmu(E_{1})+\frac{1}{2}
 \left(\frac{t'}{r'_1}-\frac{t}{r_1}\right)
 +s \leq \pmu(E_{1})+ \frac{t (r_1-r'_1)}{2r_1r'_1}+s .
 \end{eqnarray*}
Using that $\frac{t}{r_1} \leq \pmu(E_{1})-\pmu(E_{2})+2s$, by the
exact sequence (\ref{eqn:harto12}) and $1\le r'_1 \leq r_1-1$, we
get bounds on $\pmu(E'_{1})$ and $\pmu(E'_{2})$. Substituting
these bounds into (\ref{eqn:sigma}) and using that $r_1'-r_2'\geq
1$ and $r_1',r_2'\leq r_1=r_2$, we get a bound on $\sigma_c$, as
required.
\end{proof}

With this result, we may introduce the notation $\cN_L^s$ for the
moduli space of $\sigma$-stable triples for any value
$\sigma>\sigma_2$. We shall refer to this as the moduli space
\emph{for large values of $\sigma$}. There is an obvious condition
for $\cN_L^s$ to be non-empty. Let $\phi: E_2\to E_1(D)$ be a
parabolic morphism which is moreover injective. For any $p\in D$,
it induces a parabolic map $\phi_p\in \PH(E_{2,p},E_{1,p})$. This
satisfies
  $$
  \dim \im \phi_p \leq r_1-r_p,
  $$
with $r_p$ defined in (\ref{eqn:rp}) (cf.\ Remark \ref{rem:rp}).
Therefore for any parabolic map $\phi\in \PH(E_2,E_1(D))$, we have
that
 \begin{equation}\label{eqn:harto7}
 d_1+r_1 s -d_2 \geq \sum_{p\in D} r_p.
 \end{equation}
Let us now see that this is a sufficient condition for
non-emptiness and irreducibility of $\cN_L^s$. First we need some
preliminary results.

\begin{lem} \label{lem:finaldeltunel1}
 If both $E_2$ and $E_1$ are
 parabolic stable bundles, and $\phi:E_2\to
 E_1(D)$ is an injective parabolic map, then $T=(E_1,E_2,\phi)$ is
 a $\sigma$-stable triple for large values of $\sigma$.
\end{lem}

\begin{proof}
 Any subtriple $T'\subset T$ should have $r_2'\leq r_1'$.
 The stability of the bundles implies that
 $\pmu(E_1')<\pmu(E_1)$ and $\pmu(E_2')<\pmu(E_2)$, from where
 it follows that $\pmu_{\sigma}(T')<\pmu_{\sigma}(T)$, for any $\sigma$,
 in particular for large values of $\sigma$.
 \end{proof}

\begin{lem}\label{lem:harto8}
Let $L$ be a fixed parabolic line bundle. Consider the moduli
space $\cN_\sigma(r_1,r_2,d_1,d_2;\alpha,\beta)$ of
$\sigma$-stable parabolic triples $T=(E_1,E_2,\phi)$ of degrees
$(d_1,d_2)$ and weight types $(\alpha,\beta)$. Let
$(\tilde{d}_1,\tilde{d}_2)$ and $(\tilde{\alpha},\tilde{\beta})$
the degrees and weight types of the triples of the form
$(E_1\otimes^p L, E_2\otimes^p L,\phi)$. Then
$(E_1,E_2,\phi)\mapsto (E_1\otimes^p L, E_2\otimes^p L,\phi)$
gives an isomorphism $\cN_\sigma(r_1,r_2,d_1,d_2;\alpha,\beta)
\cong \cN_\sigma(r_1,r_2,\tilde{d}_1,
\tilde{d}_2;\tilde{\alpha},\tilde{\beta})$. \hfill $\Box$
\end{lem}

Let us see that tensoring with a suitable parabolic line bundle
allows us to reduce to the case $r_p=0$ for all $p\in D$. For this
we need an alternative characterization of $r_p$. Fix $p\in D$,
and denote by $\alpha_1<\cdots <\alpha_{r_1}$ the weights of
$E_{1,p}$ and by $\beta_1<\cdots<\beta_{r_1}$ the weights of
$E_{2,p}$, since $r_2=r_1$. Extend the weights to an infinite
sequence of real numbers by declaring $\alpha_{k+r_1 m}=\alpha_k +
m$, $1\leq k\leq r_1$, $m\in \ZZ$. This means that we have a
sequence
  $$
   \cdots <\alpha_{r_1}-1< \alpha_1<\cdots <
   \alpha_{r_1} <\alpha_1+1 < \alpha_2+ 1 < \cdots
  $$
In this strictly increasing sequence $\ZZ\to \RR$, $1$ is sent to
$\alpha_1$ characterized as the smallest non-negative number in
the sequence. Similarly consider the infinite sequence $\beta_k$
from the weights of $E_{2,p}$. Define the functions:
  \begin{equation} \label{eqn:harto9}
  \begin{array}{ccl}
    f: [0,\infty) & \rightarrow & \RR, \\
     x & \mapsto & \# \{\alpha_k \ | \ 0<\alpha_k< x\}, \\
    g: [0,\infty) & \rightarrow & \RR ,\\
     x & \mapsto & \# \{\beta_k \ | \ 0<\beta_k \leq x\}.
  \end{array}
  \end{equation}
Note that $f(x+1)=f(x)+r_1$ and $g(x+1)=g(x)+r_1$. Now we have

\begin{lem} \label{lem:rp-otravez}
  $r_p=\max (f-g)= \max_{[0,1)}(f-g)$.
\end{lem}

\begin{proof}
 The way $f$ and $g$ are defined, $f-g$ is a right-continuous step function, with
 jumps by $+1$ at the points $\alpha_k$
 and $-1$ at the points $\beta_k$. As
 $f-g$ is $1$-periodic, the existence of
 maximum and the equality $\max (f-g)= \max_{[0,1)}(f-g)$ are
 clear. Let $M=\max(f-g)$ and $x_0\in [0,1)$ be a point which is not
 a weight and satisfies $(f-g)(x_0)=M$. Then, writing $k=f(x_0)$,
 we have $\alpha_k <x_0< \alpha_{k+1}$
 and $k-M=g(x_0)$, i.e. $\beta_{k-M} < x_0< \beta_{k-M+1}$. The
 maximality of $f-g$ at $x_0$ implies that
 we have $\beta_{k-M} < \alpha_k <x_0< \beta_{k-M+1}<\alpha_{k+1}$.
 So any parabolic map $\phi_0:E_{2,p}\to E_{1,p}$ satisfies that
 $\phi_0(E_{2,p,k-M+1})\subset E_{1,p,k+1}$ and hence
  $$
  \dim \ker\phi_0 \geq \dim E_{2,p,k-M+1} - \dim E_{1,p,k+1}
  = (r_1-k+M)-(r_1-k)=M.
  $$

 Conversely, let $\phi_0:E_{2,p}\to E_{1,p}$ be a map such
 that $\phi_0(E_{2,p,k-M+1})\subset E_{1,p,k+1}$
 for each $k$. Then $\phi_0$ is a parabolic map: for if
 $\beta_i> \alpha_j$, take $\beta_i>x>\alpha_j$.
 So $g(x)\leq i-1$ and $f(x)\geq j$. So
 $j-i+1 \leq f(x)-g(x) \leq M$ and hence $i \geq j-M+1$. Thus
 $\phi_0(E_{2,p,i}) \subset \phi_0(E_{2,p,j-M+1}) \subset
 E_{1,p,j+1}$. On the other hand, it is clear that there are maps
 satisfying $\phi_0(E_{2,p,k-M+1})\subset E_{1,p,k+1}$
 for each $k$ with $\dim \ker \phi_0=M$. Hence there are parabolic
 maps $\phi_0$ with $\dim\ker\phi_0=M$, completing the proof that
 $M=r_p$.
\end{proof}

\begin{prop} \label{prop:harto6}
  There exists a suitable parabolic line bundle $L$ such that
  the moduli space of $\sigma$-stable triples of the form $(E_1\otimes^p L,
  E_2\otimes^p L, \phi)$ has associated
  $\tilde{r}_p=0$, for all $p\in D$.
\end{prop}

\begin{proof}
 We shall assume that there is only one point $p\in D$ and we shall tensor
 with a parabolic line bundle of the form $L=\cO_{[x]}$, i.e. the
 trivial line bundle with weight $x\in [0,1)$ at $p$.
 Take $x_0\in (0,1)$ which does not coincide with any weight and
 gives the maximum value of the function $f-g$. Let $L=\cO_{[1-x_0]}$.
 Denoting by $k_0=f(x_0)$, the weights of $E_2\otimes^p L$ are
 $$
 0\leq \alpha_{k_0+1} - x_0 < \cdots < \alpha_{r_1}-x_0 < \alpha_1
 -x_0+1 <\cdots <\alpha_{k_0} -x_0+1 <1 \,
 $$
(see \cite{ggm}). Said otherwise, if $\tilde{\alpha}_k$ is the
infinite sequence associated to the weights of
$\tilde{E}_2=E_2\otimes^p L$, then
$\tilde{\alpha}_k=\alpha_{k+k_0}-x_0$. The function $\tilde{f}$
associated to $\tilde{E}_2$ as in (\ref{eqn:harto9}) is
 \begin{eqnarray*}
 \tilde{f}(x )&=& \# \{ \tilde\alpha_k \ | \ 0<\tilde\alpha_k< x\}\\
 &=& \#\{ \alpha_k \ | \ 0 <\alpha_k -x_0< x\} \\
 &=& \#\{ \alpha_k \ | \ x_0 <\alpha_k < x+x_0\} \\
 &=& f(x+x_0)-f(x_0),
 \end{eqnarray*}
the last equality because $x_0$ is not a weight of $E_{2,p}$.
Analogously for $\tilde{E}_1=E_1\otimes^p L$, the function
$\tilde{g}$ associated to it is
 $$
 \tilde{g}(x)= g(x+x_0)-g(x_0).
 $$
Then the number $r_p$ associated to the moduli spaces of triples
$(\tilde{E}_1,\tilde{E}_2,\phi)$ is
 $$
 \tilde{r}_p=\max (\tilde{f}(x)-\tilde{g}(x))=\max (
 f(x+x_0)-g(x+x_0)) - M = 0.
 $$
\end{proof}

\begin{prop} \label{prop:harto10}
  Assume that $r_1=r_2$ and $r_p=0$ for all $p\in D$. Then the moduli space of
  $\sigma$-stable triples for $\sigma$
  large and $d_2+r_2s=d_1$ is irreducible.
\end{prop}

\begin{proof}
Any triple $T=(E_1,E_2,\phi)$ in $\cN_L^s$ satisfies that
$\phi:E_2\to E_1(D)$ is generically an isomorphism by Proposition
\ref{prop:inj}. So the condition on the degrees implies that it is
an isomorphism of bundles. Moreover, by Lemma \ref{lem:t'}, the
family $\cH$ of bundles $E_1$ appearing as part of triples of
$\cN_L^s$ is a bounded family which is irreducible and the generic
element is a stable bundle (see \cite{bgg2}).

Let us study the fibres of $\cN_L^s\to \cH$. Fix $E_1\in \cH$ and
consider the fibre over $E_1$. Identifying $E_2$ with $E_1(D)$ (as
bundles) via the isomorphism $\phi$, an element
$(E_1,E_2,\phi)=(E_1,E_1(D),\id)$ in the fibre consists on giving
for each $p\in D$ a flag for $V=E_{1,p}$ and a flag for
$V=E_{2,p}$ such that the identity map $\id: V\to V$ is a
parabolic map with respect to these flags. For simplicity, assume
there is only one point $p\in D$. Let
  $$
  \cF_1=\{0 \subset V_1 \subset V_2 \subset \cdots \subset V_{r_1}
  =E_{1,p} \  | \ \dim V_i=i \}
  $$
be the space parametrizing (complete) flags at $E_{1,p}$, with
fixed weights $\alpha_1<\cdots <\alpha_{r_1}$. This is an
irreducible variety. Analogously define the space
  $$
  \cF_2=\{0 \subset W_1 \subset W_2 \subset \cdots \subset W_{r_1}
  =E_{2,p} \ | \  \dim W_i=i \}
  $$
of (complete) flags for $E_{2,p}$, with fixed weights
$\beta_1<\cdots <\beta_{r_1}$. The condition $r_p=0$ means that
$g(x)\leq f(x)$, for all $x$, with the notation
(\ref{eqn:harto9}). The identity map is parabolic if $W_i \subset
V_{i+k(i)}$, $1\leq i \leq r_1$, for some set of integers
$k(i)\geq 0$ such that $0<1+k(1)\leq 2+k(2) \leq \cdots \leq
r_1+k(r_1)=r_1$. The set of compatible flags if given by
  \begin{equation}\label{eqn:compatibleflags}
  \cF= \{ (F_1,F_2) \ | \ W_i \subset V_{i+k(i)}, 1\leq i
  \leq r_1 \} \subset \cF_1 \times \cF_2 \, .
  \end{equation}
This is also an irreducible variety, as $\cF\to \cF_1$ is a
fibration with irreducible base and irreducible fibres. Note that
the other projection $\cF\to \cF_2$ is also surjective.

A generic stable bundle $E_1$ satisfies that a generic flag
$F_1\in\cF_1$ gives a parabolic stable bundle. Let $U_1\subset
\cF_1$ be a (dense) open subset with this property. Analogously
consider a dense open subset $U_2\subset \cF_2$ such that
$E_2=E_1(D)$ with a flag $F_2\in\cF_2$ is parabolically stable. If
$\cF\cap (U_1\times U_2)=\emptyset$ then $\cF\subset
((\cF_1-U_1)\times \cF_2 )\bigcup (\cF_1\times (\cF_2-U_2 ))$.
Being irreducible, $\cF$ should be contained in either of
$((\cF_1-U_1)\times \cF_2 )$ or $(\cF_1\times (\cF_2-U_2 ))$. This
contradicts the surjectivity of both $\cF\to \cF_1$ and $\cF\to
\cF_2$. This proves that $\cF\cap (U_1\times U_2) \neq \emptyset$,
so the generic element of $\cF$ gives parabolic stable bundles
$E_1$ and $E_2$. By Lemma \ref{lem:finaldeltunel1}, such element
is $\sigma$-stable for $\sigma$ large. Therefore the generic
stable bundle $E_1$ satisfies that the fibre of $\cN_L^s\to \cH$
is an open subset of the space of compatible flags $\cF$. This
shows that $\cN_L^s$ is irreducible and non-empty.
\end{proof}

\begin{thm} \label{thm:r1=r2}
 Suppose that $r_1=r_2$ and that $d_1+r_1s-d_2=\sum_{p\in D} r_p$.
 Then the moduli space $\cN_L^s$ is irreducible, of
 the expected dimension and non-empty.
\end{thm}

\begin{proof}
  By Proposition \ref{prop:harto6} there exists a parabolic line
  bundle $L$ such that $(E_1,E_2,\phi)\mapsto (\tilde{E}_1=E_1\otimes^p L,
  \tilde{E}_2=E_2\otimes^p
  L,\phi)$ gives an isomorphism of moduli spaces of $\sigma$-stable
  triples $\cN_\sigma(r_1,r_1,d_1,d_2;\alpha,\beta) \cong
  \cN_\sigma(r_1,r_1,\tilde{d}_1,\tilde{d}_2;\tilde\alpha,\tilde\beta)$ such that
  $\tilde{r}_p=0$ for each $p\in D$. Then
  $$
  \tilde{d}_1+r_1s -\tilde{d}_2 =d_1+r_1s-d_2 - \sum_{p\in D} r_p.
  $$
  This is easily seen by computing the degrees $\tilde{d}_1$ and
  $\tilde{d}_2$. For instance, suppose that there is only one
  point $p\in D$. Then, with the notations of the proof of
  Proposition \ref{prop:harto6},
   \begin{eqnarray*}
   \tilde{d}_1 &=& \deg \tilde{E}_1=\pdeg (E_1\otimes^p L)- \sum
   \tilde{\alpha}_k \\
   &=& \pdeg (E_1) + r_1 \pdeg (L) - \left(\sum (\alpha_k - x_0)
   +k_0\right) \\
   &=& d_1 +\sum \alpha_k + r_1 (1-x_0) - \sum \alpha_k + r_1
   x_0-k_0 \\
   &=& d_1+r_1 - k_0=d_1+r_1-f(x_0)\, .
  \end{eqnarray*}
Analogously, $\tilde{d}_2=d_2+r_1-g(x_0)$, so that $\tilde{d}_1
-\tilde{d}_2 =d_1-d_2 - r_p$.

Now the moduli space $\cN_L^s(r_1,r_1,\tilde{d}_1,\tilde{d}_2;
\tilde{\alpha},\tilde{\beta})$ is non-empty and irreducible by
Proposition \ref{prop:harto10}. So the same is true of our initial
moduli space by using Lemma \ref{lem:harto8}. The dimension
statement follows from Corollary \ref{cor:expected}.
\end{proof}

\begin{thm} \label{thm:harto11}
Suppose that $r_{1}=r_{2}$ and $d_{1}+r_{1}s-d_{2}\ge \sum_{p\in
D} r_{p}$. Then the moduli space $\cN_{L}^s$ is non-empty, of the
expected dimension and irreducible.
\end{thm}

\begin{proof}
The dimension statement follows from Corollary \ref{cor:expected}.
Arguing as in the proof of Theorem \ref{thm:r1=r2}, we may suppose
that $r_{p}=0$, for $p\in D$. Now, there exist triples
$\phi:E_{2}\to E_{1}(D)$, with $\phi$ injective, $E_{1}$ and
$E_{2}$ stable bundles, and satisfying that the torsion sheaf
quotient of the map $\phi$ is generic (in particular, supported on
$X \setminus D$). This follows from \cite{bgg2}, where non
parabolic $\sigma$-stable triples for $\sigma$ large are found by
constructing $\sigma$-stable triples with these properties.

Now the argument of the proof of Proposition \ref{prop:harto10}
works here to find parabolic structures on $E_{1}$ and $E_{2}$
such that $(E_{1},E_{2},\phi)$ is a $\sigma$-stable parabolic
triple for $\sigma$ large, since the only necessary fact is that
$\phi_{p}:E_{2,p}\to E_{1,p}$ is an isomorphism for all $p\in D$.
This gives the non-emptiness of $\cN_L^s$.

For proving the irreducibility of $\cN_L^s$, the main obstacle are
the triples with quotient supported at points of $D$. We work as
follows. Let $\cH$ be the family of bundles $E_{1}$ appearing in
triples $T=(E_1,E_2,\phi) \in \cN_{L}^s$. This is a bounded and
irreducible family whose generic element $E_{1}\in \cH$ is a
generic stable bundle. Let $\cQ=\Quot^t(\cH)$ be the Quot scheme
parametrizing quotients $E_1(D) \to S$, with $E_1\in \cH$ and
$t=\length S=d_{1}+r_{1}s-d_{2}-\sum_{p\in D} r_{p}$. The kernel
of a generic element in $\cQ$ is a stable bundle $E_2$. If the
support of $S$ is contained in $X\setminus D$, then the fiber of
the map $\cN_L^s\to \cQ$ over a quotient $E_1(D)\to S$ in $\cQ$ is
a subset of the set of compatible flags $\cF$ defined in
(\ref{eqn:compatibleflags}). For a generic element in $\cQ$, this
is actually an open subset of $\cF$, as proved in the proof of
Proposition \ref{prop:harto10}. This produces an open subset
$U\subset \cN_L^s$, which is of dimension
  $$
  \dim \cQ + \dim \cF\, .
  $$

Let us see the irreducibility of $\cN_L^s$ by checking that $\dim
(\cN_L^s\setminus U)<\dim U$. Certainly, the only effect that we
must take care of is the jumping in the dimension of the fiber of
$\cN_L^s\to \cQ$ when the torsion sheaf is supported at some
points of $D$. Let $p\in D$, and suppose that $p$ is in the
support of $S$, say $S_{p}=\CC^{l}$. The set of quotients
$E_{1,p}\to S_{p}$ is parametrized by the grassmannian
$\Gr(l,r_{1})$. The codimension of the space $\cQ^l\subset \cQ$
parametrizing such quotients is
 $$
 r_{1}\length S-(r_{1}(\length
 S-l)+l(r_{1}-l))=r_{1}l-lr_{1}+l^{2}=l^{2}.
 $$
Now let us compute the dimension of the fiber of $\cN_L^s\to \cQ$
over a point in $\cQ^l$. With the definition of $k(i)$ given in
Proposition \ref{prop:harto10}, such fiber is the space
 $$
 \cF_{\ast} = \{(W_{i},V_{i})\in \cF_{1}\times \cF_{2}
 \ |\ \phi(W_{i})\subset V_{i+k(i)}\} \, .
 $$
Equivalently, $(W_{i},V_{i})\in \cF_{\ast}\Leftrightarrow
W_{i}\subset\phi^{-1}(V_{i+k(i)})$. It remains to see that
  $$
  \dim \cF_{\ast}-\dim \cF<l^{2}\, .
  $$

The fibration $\cF\to\cF_{1}$ is surjective and the dimension of
the fiber is
 $$
 \sum_{i=1}^{r_1}k(i)
 $$
Let us compute the dimension of a fiber of $\cF_{\ast}\to \cF_1$.
Such dimension depends on the flag $\{V_{i}\}\in \cF_{1}$, so we
need to stratify $\cF_1$ as follows. The flag $\{V_{i}\}$ is
determined by a collection of numbers $0\le a_{1}\le \ldots\le
a_{r_1}=r_1-l$ such that
  $$
  \begin{array}{ccccccc}
  0&\subset & V_{1}\cap \img(\phi)&\subset &\cdots &\subset &V_{r_1}\cap
  \img(\phi)=\img(\phi) \\
   \| &&\| && &&\| \\
  0&\subset &\CC^{a_{1}}&\subset &\cdots&\subset&\CC^{a_{r_1}}=\CC^{r_1-l}
  \end{array}
  $$
Clearly, $a_{i+1}=a_{i}+\delta_{i+1}$ ($a_0=0$) where there are
uniquely defined $1\leq i_{1}<\ldots <i_{r_1-l}\leq r_1$ such that
$\delta_{i_{k}}=1$ and $\delta_j=0$ for $j\neq i_k$, $k=1,\ldots,
r_1-l$. The codimension of the stratum $S_{a_1,\ldots,a_{r_1}}
\subset \cF_1$ defined by such $\{V_{i}\}$ is
  $$
  \sum_{k=1}^{r_1-l}(l-i_{k}+k).
  $$

The fiber of $\cF_{\ast}\to \cF_{1}$ over $\{V_{i}\}\in
S_{a_1,\ldots,a_{r_1}}$ is given by flags $\{W_{i}\}\in \cF_2$
such that $W_{i}\subset \tilde{V}_{i+k(i)}$, with
$\tilde{V}_{i}=\phi^{-1}(V_{i})\cong \CC^{l+a_{i}}$. The dimension
of such fiber is thus
  \begin{eqnarray*}
  \sum_{i=1}^{r_1}(l+a_{i+k(i)}-i)&\le& \sum_{i=1}^{r_1}(l+a_{i}-i)+\sum
  k(i)\\
  &=&\sum_{i=1}^{r_1}(l-i)+\sum_{k=1}^{r_1-l}(r_1-i_{k}+1)+\sum k(i)
  \end{eqnarray*}
So the dimension of the preimage of $S_{a_1,\ldots,a_{r_1}}$ by
the map $\cF_{\ast}\to\cF_1$ is less than or equal to
  \begin{align*}
  \dim
  \cF_{1} &-\sum_{k=1}^{r_1-l}(l-i_{k}+k)+ \sum_{i=1}^{r_1}(l-i)+
  \sum_{k=1}^{r_1-l}(r_1-{i_{k}}+1)+\sum k(i) \notag \\
  &=\dim \cF_{1}+\sum k({i}) +\frac{l^{2}-l}{2}= \dim
  \cF +\frac{l^{2}-l}{2}.
  \end{align*}
Since this is true for any stratum, we have
 $$
 \dim \cF_{\ast}\le \dim \cF +\frac{l^{2}-l}{2}<\dim \cF+ l^{2},
 $$
as required.
\end{proof}

Combining Theorem \ref{thm:harto11} with Proposition
\ref{prop:equalrs} we have the following.

\begin{cor}\label{cor:irreducibility}
Let $g>0$, $r_1=r_2$ and $d_1+r_1s-d_2 \geq \sum_{p\in D} r_p$.
Then the moduli spaces $\cN_\sigma$ are non-empty, irreducible and
of the expected dimension for any  $\sigma\geq 2g-2$.
\end{cor}

\begin{rmk} \label{rem:irreducibility}
Corollary \ref{cor:irreducibility} and the correspondence in
Proposition \ref{prop:iden} gives that the moduli space
$\cU(p,p,a,b;\alpha,\beta)$ is non-empty and connected if and only
if the following is satisfied:
 \begin{itemize}
 \item[(i)] In the case $\tau<0$. It must be $|\tau|\leq \tau_M$
 by Proposition \ref{prop:harto11}. Also, defining
 $r_x=\min \{\dim \coker\phi \ |
 \ \phi\in \PH(V_x,W_x)\}$, for $x\in D$, we must have $b+(2g-2+s)p -a\geq
 \sum_{x\in D} r_x$, by Corollary \ref{cor:irreducibility}. But
 this last condition is redundant: $\tau<0$ is equivalent to
 $\pmu (V) < \pmu(W)$, hence
  $$
  a=\deg (V) \leq \pdeg(V)<\pdeg(W)<\deg (W)+ps= b+ps +(2g-2)s,
  $$
 since $g>0$. Also, we may tensor with a suitable parabolic line
 bundle $L$ to arrange $r_x=0$, for all $x\in D$, by Proposition
 \ref{prop:harto6} (this does not change $\tau$ or the inequality
 that we need to check). So $b+(2g-2+s)p-a\geq 0$, as required.
 \item[(ii)] The case $\tau>0$ is worked out similarly, and the only
 condition we obtain is $|\tau|\leq \tau_M$.
 \end{itemize}
Note that the genericity of the weights (Assumption
\ref{assumption}) prevents the case $|\tau|=\tau_M$ to happen.
\end{rmk}

\section{Representations of  fundamental groups in $\U(p,q)$}
\label{representations}

Let $X$ be a compact Riemann surface  of genus $g\ge 0$ and let
$S=\{x_{1}, \ldots, x_{s}\}$ be a set of distinct points of $X$.
Let $\Gamma=\pi_1(X\setminus S)$ be the fundamental group of
$X\setminus S$. The group $\Gamma$ is generated by the usual
generators $a_i,b_i$, $1\leq i\leq g$, of $\pi_1(X)$, together
with additional generators  $\gamma_1,\ldots, \gamma_s$
corresponding to loops enclosing each $x_i$ simply, not enclosing
any $x_j$, $j\neq i$, and which are homotopic to zero relatively
to the base point on $X$. There is also the relation
$[a_1,b_1]\cdots [a_g,b_g]\gamma_1\cdots \gamma_s=1$, where
$[a_i,b_i]$ is the commutator of $a_i$ and $b_i$.

Parabolic Higgs bundles are related to representations of $\Gamma$.
To be precise, let us fix integers $n=\rk E$, $d=\deg E$ and
the weight type $\alpha=\{\alpha(x)\}_{x\in S}$, where
$\alpha(x)=(\alpha_{1}(x),\ldots , \alpha_{r(x)}(x))$ are weights
with multiplicities $k_i(x)$ for every $x\in S$.
It is  convenient to repeat each weight according
to its multiplicity, by setting
$\tilde \alpha_1(x)=\ldots =\tilde\alpha_{k_1(x)}(x)=\alpha_1(x)$, etc.,
thus having  weights
$0\leq \tilde \alpha_{1}(x)\le \ldots \le\tilde \alpha_n(x) <1$
(see Section \ref{sec:intro}).

For every $x_i\in S$ there is a $C_i\in \U(n)$ defined by
\begin{equation}\label{holonomy}
  C_i =
  \begin{pmatrix}
    \exp(2\pi\sqrt{-1}\tilde{\alpha}_1(x_i)) & & 0 \\
     & \ddots & \\
    0 & & \exp(2\pi \sqrt{-1}\tilde{\alpha}_n(x_i)) \\
  \end{pmatrix}.
\end{equation}
Consider the set of representations
$\Hom^{+}_{\alpha}(\Gamma,\GL(n,\CC))$ defined by semisimple
homomorphisms $\rho:\Gamma\rightarrow \GL(n,\CC)$ such that
$\rho(\gamma_i)$ is conjugated to $C_i$ by an element in
$\GL(n,\CC)$ for $1\leq i\leq s$. Here by semisimple we mean that
$\rho$ is a direct sum of irreducible representations. The moduli
space of representations of $\Gamma$ in $\GL(n,\CC)$ with fixed
holonomy in the conjugacy class of $C_i$,  is defined by the
quotient
 $$
 \cR(n;\alpha):=\frac{\Hom^{+}_\alpha(\Gamma,\GL(n,\CC))}{\GL(n,\CC)},
 $$
where $\GL(n,\CC)$ acts by conjugation. The set $\cR(n;\alpha)$
has a natural structure of a complex algebraic variety. The
following is proved by Simpson in \cite{s2}.

\begin{thm}\label{nahodge-gln}
Let $(n,d;\alpha)$  be such that
 $$
 d+ \sum_{x\in S} (\tilde \alpha_1(x)+\ldots + \tilde \alpha_n(x))=0,
 $$
i.e., the parabolic degree vanishes. Then there is a
homeomorphism
 $$
 \cR(n;\alpha) \cong  \cM(n,d;\alpha).
 $$
\end{thm}

This generalizes the theorem of Metha--Seshadri \cite{ms} which
identifies the moduli space of parabolic bundles of type
$(n,d,\alpha)$ with vanishing parabolic degree with the moduli
space of representations of $\Gamma$ in $\U(n)$ with fixed
holonomy conjugated to $C_i$ around the marked points.

There is a similar correspondence between representations of $\Gamma$ in
$\U(p,q)$ and parabolic $\U(p,q)$-Higgs bundles. To explain this,
let us come back to the notation in
Section \ref{sec:PHB} and  fix the types of the parabolic bundles $V$ and $W$
to be  $(p,a,\alpha)$ and $(q,b,\alpha')$, respectively.
For every $x_i\in S$ there are  matrices   $C_i\in \U(p)$ and
$C_i'\in \U(q)$  defined as in (\ref{holonomy}) by  the weight systems
$\alpha$ and $\alpha'$, respectively.

Consider now the set of representations
$\Hom^{+}_{\alpha,\alpha'}(\Gamma,\U(p,q))$ defined by semisimple
homomorphisms $\rho:\Gamma\rightarrow \U(p,q)$ such that
$\rho(\gamma_i)$ is conjugated to $C_i\times C_i'\in \U(p)\times
\U(q)$ (recall that $\U(p)\times \U(q)$ is the maximal compact
subgroup of $\U(p,q)$) by an element in $\U(p,q)$ for $1\leq i\leq
s$. Define the moduli space of representations of $\Gamma$ in
$\U(p,q)$ with fixed holonomy $\U(p,q)$-conjugated to  $C_i\times
C_i'$  by the quotient
 $$
 \cR(p,q;\alpha,\alpha')
 :=\frac{\Hom^{+}_{\alpha,\alpha'}(\Gamma,\U(p,q))}{\U(p,q)}.
 $$

The set $\cR(p,q;\alpha,\alpha')$ is a real analytic variety.
We can adapt the arguments of Simpson \cite{s2} to prove  the following.
\begin{thm}\label{nahodge-upq}
Let $(p,a,\alpha)$ and $(q,b,\alpha')$ be such that
 $$
 \pdeg(V)+\pdeg(W)=a+b+ \sum_{x\in S} (\tilde \alpha_1(x)
 +\ldots + \tilde \alpha_p(x) + \tilde \alpha'_1(x)  +\ldots +  \tilde \alpha'_q(x))=0.
 $$
Then there is a homeomorphism
 $$
 \cR(p,q;\alpha,\alpha')\cong \bigsqcup_{a,b} \cU(p,q,a,b;\alpha,\alpha').
 $$
\end{thm}

Note that $(p,q,a,b;\alpha,\alpha')$ must also  satisfy the
Milnor--Wood inequality, which in these cases  reduces to
  $$
  |\pdeg(V)|\leq \min\{p,q\}(g-1+s/2),
  $$
since $\pdeg(W)=-\pdeg(V)$.

Combining Theorem \ref{nahodge-upq} and Theorem
\ref{thm:connected-minima} we have the following.

\begin{thm}\label{nahodge-upq2}
Under the genericity conditions given by Assumption
\ref{assumption}, and for $g>0$, the number of non-empty connected
components of $\cR(p,q;\alpha,\alpha')$ equals the number of
integers $a$ such that
 $$
 |a+ \sum_{x\in S} (\tilde \alpha_1(x)+\ldots + \tilde
 \alpha_p(x))| \leq\tau_L/2,
 $$
where $\tau_L$ is given by (\ref{eqn:tauL}).
\end{thm}

\begin{rmk}
 The condition on the genus $g$ comes from Theorem \ref{nahodge-upq}.
\end{rmk}

Like in the proof of Theorem \ref{nahodge-gln} (\cite{s2}), the
main ingredients in the proof of Theorem \ref{nahodge-upq} are, on
the one hand, the correspondence given by Theorem \ref{hk} between
polystable parabolic $\U(p,q)$-Higgs bundles and solutions to
Hitchin equations, and, on the other, the existence of a harmonic
adapted metric on a $\U(p,q)$-bundle with a semisimple meromorphic
flat connection with simple poles.   To see  this, let us come
back to the framework of Section \ref{gauge}, and consider smooth
parabolic vector bundles $V$ and $W$ of types $(p,a;\alpha)$ and
$(q,b;\alpha')$, respectively. On the bundle $V\oplus W$ we
consider flat $\U(p,q)$-connections $D$ on $X\setminus S$,
meromorphic at $x_i\in S$ and whose residue at $x_i$ is conjugated
to $C_i\times C_i'$. We say that $D$ is semisimple if the
corresponding representation is semisimple. These  connections are
in correspondence with elements in
$\Hom^{+}_{\alpha,\alpha'}(\Gamma,\U(p,q))$.

Let $h=(h_V,h_W)$, where $h_V$ and $h_W$ are adapted hermitian
metrics on $V$ and $W$, respectively. We decompose $D$ as
$D=d_A+\Psi$, where $d_A$ is a $\U(p)\times \U(q)$ connection  and
$\Psi$ takes values in $\lie{m}$, where
$\lie{u}(p,q)=\lie{u}(p)\oplus \lie{u}(q) + \lie{m}$ is the Cartan
decomposition of the Lie algebra of $\U(p,q)$.  We say that $h$ is
harmonic if $d_A^\ast \Psi=0$. Then the following can be proved
easily adapting the results in \cite{c,s2}.

\begin{thm}\label{harmonic}
A connection $D$ as above  is semisimple if and only if there exists a harmonic
hermitian metric $h=(h_V,h_W)$.
\end{thm}
The relation with parabolic $\U(p,q)$-Higgs bundle is given as follows.
If $D$ is semisimple flat connection as above   and $h$ is a harmonic
solution, then the pair $(d_A,\Phi)$, where $\Phi$ is determined
 by the equation $\Psi=\Phi +\Phi^*$,
solves the $\U(p,q)$-Hitchin equations  and hence, by Theorem \ref{hk},
corresponds to a polystable parabolic $\U(p,q)$-Higgs bundle. Conversely,
if we have a polystable parabolic $\U(p,q)$-Higgs bundle we can find a solution
$(d_A,\Phi)$ to the Hitchin equations, and then out of it a solution
to the harmonic equation on the flat connection $D=d_A +\Phi+\Phi^*$,
which is then semisimple by Theorem \ref{harmonic}.

\section{Elliptic surfaces, orbifolds and parabolic Higgs bundles}

Parabolic bundles have been related by several authors to unitary
representations of the fundamental group of elliptic surfaces of
general type (\cite{Ba,SS}). The key fact is that the fundamental
group of such a surface is isomorphic to the orbifold fundamental
group of an orbifold Riemann surface, whose unitary
representations are, in turn,  related to parabolic bundles by the
Metha--Seshadri theorem \cite{ms,b,bo,nst2}.

Let $X$ be a compact Riemann surface  of genus $g\ge 0$ and let
$S=\{x_{1}, \ldots, x_{s}\}$ be a set of distinct points of $X$.
Suppose that for each $i$ we are given integers $m_i\geq 1$, such
that $2g+\sum_{1\leq i\leq s}(1-1/m_i)>2$. We call the data of
$X$, $S$, and $m_i$,  $1\leq i\leq g$, a $2$-{\em orbifold}. As in
Section \ref{representations}, let $\Gamma=\pi_1(X\setminus S)$ be
the fundamental group of $X\setminus S$. As we have seen in
Section \ref{representations}, $\Gamma$ has $2g+s$ generators
$a_i,b_i$, $1\leq i\leq g$, and $\gamma_j$, $1\leq j\leq s$,
satisfying the relation
 $$
 \prod_{1\leq i\leq g}[a_i,b_i]\cdot\prod_{1\leq j\leq s}\gamma_j=1.
 $$
We define the {\em orbifold fundamental group} $\pi_1^{\orb}(X)$
as  the quotient of $\Gamma$ by the smallest normal subgroup
containing  $\gamma_i^{m_i}$. Thus $\pi_1^{\orb}(X)$ is freely
generated by the elements $a_i,b_i$, $1\leq i\leq g$, and
$\gamma_j$, $1\leq j\leq s$, subject to the relations
 $$
 \prod_{1\leq i\leq g}[a_i,b_i]\cdot \prod_{1\leq j\leq s}\gamma_j=1,\;\;\;
 \mbox{and}\;\;\; \gamma_j^{m_j}=1, \;\;\; 1\leq j\leq s.
 $$
The $2$-orbifold Riemann surface ought to be thought of  as
a Riemann surface with singularities at
the points $x_i$,  which locally are  of the form $\Delta/\ZZ_{m_i}$,
where $\Delta$ is the unit disc in $\CC$.
The group $\pi_1^{\orb}(X)$ is clearly the fundamental group of this
orbifold surface (see \cite{bo,nst} and references there for basic
facts on orbifold surfaces).

The following is proved  in \cite{Dol,Ue} (see also \cite{Fri,SS}).
\begin{thm}
Given an orbifold fundemental group $\pi_1^{\orb}(X)$
and an integer $\chi>0$, there is
an elliptic surface $Y$, unique up to diffeomorphism,  with
 $$
 \pi_1(Y)=\pi_1^{\orb}(X),\;\;\;\mbox{and}\;\;\;\chi(\cO_Y)=\chi.
 $$
Conversely,  given an elliptic surface $Y$ with $b_1(Y)$ even,
$\chi(\cO_Y)>0$ and $\kod(Y)=1$ we have
 $$
 \pi_1(Y)=\pi_1^{\orb}(X),
 $$
for some $2$-orbifold Riemann surface $X$.
\end{thm}

To understand this result and the relation of $Y$ to the
$2$-orbifold $X$, recall that an elliptic surface is a smooth
compact complex surface $Y$ with a fibration $f:Y\rightarrow X$
onto a Riemann surface $X$ such that the generic fibre is an
elliptic curve (the complex structure of the fibre may vary from
point to point). In some special points the fibre may degenerate
into nodal fibers. This is always the case for the elliptic
surfaces we are dealing with. Technically this is the condition
$\chi>0$. The effect of these singularities is that they kill the
extra generators of the fundamental group determined by the fibre.
In addition  to these nodal fibres there are multiple fibres,
located over the marked points of $X$. They are defined
analogously to orbifold singularities: a neighbourhood $Y_m$ of
such a multiple fibre in $X$ is the quotient by a finite cyclic
group,
 $$
 f: Y_m\cong (\Delta \times E_{\tau(z)})/\ZZ_m \lto
 \Delta/\ZZ_m\cong \Delta
 $$
defined by $[(t,c)]\mapsto t^m=z$, where $\Delta$ is the unit disc
in $\CC$, $E_\tau$ is the torus $\CC/\ZZ \oplus \ZZ \tau$, and the
generator of $\ZZ_m$ acts as $(t,c)\mapsto (t\cdot
\exp(2\pi\sqrt{-1}/m), c+1/m)$. The crutial difference of a
multiple fibre of $Y$ and the orbifold point is, however, that
this action is free and hence the quotient is smooth. Roughly
speaking, the orbifold singularity is now hidden in the map $f$
between two smooth manifolds $Y$ and $X$.

To relate representations  $\rho: \pi_1^{\orb}(X)\to \GL(n,\CC)$
to parabolic Higgs  bundles, we observe  that $\rho(\gamma_i)$
must be conjugated to a matrix of the form
\begin{equation}\label{holonomy2}
  C_i =
  \begin{pmatrix}
    \exp(2\pi\sqrt{-1}\frac{l_1(x_i)}{m_i}) & & 0 \\
     & \ddots & \\
    0 & & \exp(2\pi \sqrt{-1}\frac{l_n(x_i)}{m_i}) \\
  \end{pmatrix}
\end{equation}
for integers $l_j(x_i)$ such that
 \begin{equation}\label{integers}
 0\leq l_1(x_i)\leq \ldots \leq l_n(x_i)<m_i.
 \end{equation}
This follows from the fact that  ${\rho(\gamma_i)}^{m_i}=I$. Such
a representation of $\pi_1^{\orb}(X)$ lifts to a representation
$\tilde \rho: \Gamma \to  \GL(n,\CC)$. Conversely, if $\tilde
\rho:\Gamma \to  \GL(n,\CC)$ is such that $\rho(\gamma_i)$ is
conjugated to a matrix $C_i$ as above then $\tilde\rho$ descends
to a representation $\rho: \pi_1(X^{\orb})\to \GL(n,\CC)$. We thus
have proved the following.

\begin{prop}\label{orbi-gln}
There is a one-to-one correspondence between representations
$\rho: \pi_1(X^{\orb})\to \GL(n,\CC)$ and representations
$\tilde \rho: \Gamma \to  \GL(n,\CC)$ such that $\tilde \rho (\gamma_i)$
is conjugated to a matrix of the form (\ref{holonomy2})
for integers  $l_j(x_i)$ satisfying (\ref{integers}).
\end{prop}

Similarly, we have the following.

\begin{prop}\label{orbi-upq}
There is a one-to-one correspondence between representations
$\rho: \pi_1^{\orb}(X)\to \U(p,q)$ and representations $\tilde
\rho :\Gamma \to  \U(p,q)$ such that $\tilde \rho (\gamma_i)$ is
$\U(p,q)$-conjugated to an element of the form $C_i\times
C_i'\subset \U(p)\times \U(q)$ with $C_i$ and $C_i'$ like in
(\ref{holonomy2}), defined for integers $l_j(x_i)$ and $l'_k(x_i)$
satisfying
\begin{equation}\label{integers2}
 0\leq l_1(x_i)\leq \ldots \leq l_p(x_i)<m_i \;\;\;  \mbox{and}\;\;\;
 0\leq l'_1(x_i)\leq \ldots \leq l'_q(x_i)<m_i.
\end{equation}
\end{prop}

Let
 \begin{equation}\label{lambda}
 \lambda =\{\lambda(x_i)=(l_1(x_i),\ldots, l_n(x_i))\}_{x_i\in S},
 \end{equation}
where $l_j(x_i)$ are  integers  satisfying (\ref{integers}). Let
$\cR^{\orb}_X(n;\lambda)$ and $\cR_Y(n;\lambda)$ be  the moduli
spaces of semisimple representations of $\pi_1^{\orb}(X)$ and
$\pi_1(Y)$ in $\GL(n,\CC)$ such that $\rho(\gamma_i)$ is
conjugated to the matrix (\ref{holonomy2}). Similarly, let
 \begin{equation}\label{lambda2}
 \lambda =\{\lambda(x_i)=(l_1(x_i),\ldots, l_p(x_i))\}_{x_i\in S}
 \;\;\; \mbox{and}\;\;\;
 \lambda' =\{\lambda'(x_i)=(l_1'(x_i),\ldots, l'_q(x_i))\}_{x_i\in S}
 \end{equation}
satisfying (\ref{integers2}). Let
$\cR^{\orb}_X(p,q;\lambda,\lambda')$ and
$\cR_Y(p,q;\lambda,\lambda')$ be the moduli spaces of semisimple
representations of $\pi_1^{\orb}(X)$ and $\pi_1(Y)$ in $\U(p,q)$
such that $\rho(\gamma_i)$ is conjugated to a matrix $C_i\times
C_i'$ like in Proposition \ref{orbi-upq}. Of course, since
$\pi_{1}^{\orb}(X)\cong \pi_1(Y)$,
$\cR^{\orb}_X(n;\lambda)\cong\cR_Y(n;\lambda)$ and
$\cR^{\orb}_X(p,q;\lambda,\lambda')\cong\cR_Y(p,q;\lambda,\lambda')$.

Combining Propositions
 \ref{orbi-gln} and \ref{orbi-upq} and
Theorems \ref{nahodge-gln} and \ref{nahodge-upq} we have the following.

\begin{thm} \label{elliptic-surface}
Let $\lambda$ given by (\ref{lambda}) satisfying
(\ref{integers}) and let
$\tilde \alpha(x_i) =\lambda(x_i)/m_i$. Let $(n,d)$ be such that
 $$
 d+ \sum_{x\in S} (\tilde \alpha_1(x)+\ldots + \tilde \alpha_n(x))=0.
 $$
Then
 $$
 \cR^{\orb}_X(n;\lambda)\cong\cR_Y(n;\lambda)\cong \cR(n,d;\alpha)
 \cong \cM(n,d;\alpha).
 $$
Similarly, let $\lambda$ and $\lambda'$ given by (\ref{lambda2})
satisfying (\ref{integers2}) and let $\tilde \alpha(x_i)
=\lambda(x_i)/m_i$ and $\tilde \alpha'(x_i) =\lambda'(x_i)/m_i$.
Let $(p,q,a,b)$ be such that
 $$
 a+b+ \sum_{x\in S} (\tilde \alpha_1(x) +\ldots + \tilde \alpha_p(x) +
 \tilde \alpha'_1(x)  +\ldots +  \tilde \alpha_q'(x))=0.
 $$
Then
 $$
 \cR^{\orb}_X(p,q;\lambda,\lambda')\cong
 \cR_Y(p,q;\lambda,\lambda')\cong \cR(p,q;\alpha,\alpha')\cong
 \bigsqcup_{a,b} \cU(p,q,a,b;\alpha,\alpha').
 $$
\end{thm}

As established by Simpson and Corlette, higher dimensional
non-abelian Hodge theory  (\cite{s1,c}) gives a correspondence
between semisimple  flat bundles or representations of the
fundamental group of a compact K\"ahler manifold $(Y,\omega)$, and
polystable Higgs bundles on $(Y,\omega)$ with vanishing first and
second Chern classes (see \cite{s1} for the definition of
stability). Now, a $\GL(n,\CC)$-Higgs bundle on $Y$ is defined as
a pair $(E,\Phi)$ consisting of a  holomorphic vector bundle $E$
over $Y$ and a homomorphism $\Phi: E \to E\otimes \Omega^1_Y$ such
that $[\Phi,\Phi]=0$, where $\Omega^1_Y$ is the bundle of
holomorphic  one-forms on $Y$. If $E=V\oplus W$, where $V$ and $W$
are holomorphic bundles of ranks $p$ and $q$ respectively,  and
  $$
  \Phi=\left(\begin{array}{ll}  0 & \beta \\ \gamma & 0
  \end{array}\right) :(V\oplus W)\to (V\oplus W)\otimes \Omega^1_Y,
  $$
then $(E,\Phi)$ is said to be a $\U(p,q)$-Higgs bundle. Of course,
when $Y$ is a Riemann surface we recover the original definition
of Higgs bundle since $\Omega^1_Y$ is the canonical bundle and the
condition $[\Phi,\Phi]=0$ is trivially satisifed.

If $Y$ is a complex elliptic surface as above, equipped with a
K\"ahler metric  $\omega$, non-abelian Hodge theory on $(Y,\omega)$ combined  with
Theorem \ref{elliptic-surface} gives the following.

\begin{thm}
There is a one-to-one correspondence between the moduli space of
polystable $\GL(n,\CC)$-Higgs bundles on $(Y,\omega)$ with vanishing Chern classes
and the moduli space
of parabolic $\GL(n,\CC)$-Higgs bundles on $X$ with parabolic structure on
the orbifold points.

Similarly, there is a one-to-one correspondence between the moduli space of
polystable $\U(p,q)$-Higgs bundles on $(Y,\omega)$ with vanishing Chern
classes  and the moduli space
of parabolic $\U(p,q)$-Higgs bundles on $X$ with parabolic structure on
the orbifold points.

\end{thm}

\begin{rmk}
It would be very interesting to work out this correspondence directly in a
similar fashion to what is done by Bauer \cite{Ba} for the case of moduli
spaces of vector bundles. We plan to come back to this problem in a future
paper.
\end{rmk}

\end{document}